\documentclass[twocolumn,a4paper]{article} 

\usepackage{amssymb,latexsym}
\usepackage{xcolor}
\RequirePackage[OT1]{fontenc}
\usepackage{amsthm,amsmath,natbib,xfrac}

\usepackage{bbm}
\usepackage{comment}

\RequirePackage[colorlinks,citecolor=blue,urlcolor=blue]{hyperref}

\numberwithin{equation}{section}
\theoremstyle{plain}

\newtheorem{es}{Example}[section]

\newtheorem{demo}{Proof sketch}[section]

\begin{document}

\title{Likelihood Asymptotics in Nonregular Settings: A Review with Emphasis on the Likelihood Ratio}

\author{Alessandra R. Brazzale\footnote{Alessandra R. Brazzale is Associate Professor of Statistics, Department of Statistical Sciences, University of Padova, Via Cesare Battisti 241, 35121 Padova, Italy, {\tt alessandra.brazzale@unipd.it}} ~and Valentina Mameli\footnote{Valentina Mameli is Associate Professor at the Department of Economics and Statistics, University of Udine, Via Tomadini 30/A, 33100 Udine, Italy, {\tt valentina.mameli@uniud.it}.}}

\date{April 19, 2023 \\[5ex]}

\maketitle

\begin{abstract} 
This paper reviews the most common situations where one or more regularity conditions which underlie classical likelihood-based parametric inference fail.  We identify three main classes of problems: boundary problems, indeterminate parameter problems---which include non-identifiable parameters and singular information matrices---and change-point problems.  The review focuses on the large-sample properties of the likelihood ratio statistic.  We emphasize analytical solutions and acknowledge software implementations where available.  We furthermore give summary insight about the possible tools to derivate the key results.   Other approaches to hypothesis testing and connections to estimation are listed in the annotated bibliography of the Supplementary Material. \\[1ex]
\end{abstract}

\section{Introduction} 
The likelihood ratio or Wilks statistic is the oldest of the three classical approaches of likelihood-based inference to hypothesis testing, which include the asymptotically equivalent Wald and score statistics.  Modern quests still advocate Wilks’ test to identify new rare events and/or to detect a sudden change in a data generation process.  It is commonly believed that under the null hypothesis its finite-sample distribution approaches the chi-squared distribution as the sample size goes to infinity.  However, in order to hold true, the application of Wilks’ (1938) \nocite{Wilks38} theorem to calculate measures of significance, and of subsequent and related results of large sample asymptotic theory, requires a number of regularity conditions which are often not met.  Though there has always been awareness of nonregular problems, especially practitioners may be less familiar with the resulting limiting distributions and how these are derived.  The aim of this paper is to resume a seemingly bygone statistical problem, whose inferential issues are highly relevant for a number of modern applications; see, for instance, \cite{Algeri++20}.  We present an overview of the situations in which Wilks’ (1938) theorem can fail and of how to construct valid inferences.  We will focus on the likelihood ratio statistic and its limiting distribution motivated by its widespread use for hypothesis testing, model selection and other related uses.  Analogies with alternative test statistics and/or nonparametric and semiparametric models within the frequentist but also Bayesian paradigm are listed in the annotated bibliography of the Supplementary Material.

Asymptotic theory is an essential part of statistical methodology.  It provides first thing approximate answers where exact ones are unavailable.  Beyond this, it serves to check if a proposed inferential solution provides a sensible answer when the amount of information in the data increases without limit.  Given the tremendous advances in computer age statistical inference \citep{EfronHastie16} one could be tempted to by-pass the often rather demanding algebraic derivations of asymptotic approximation.  Gaining insight in what happens to the limiting distribution of likelihood-based test statistics when one or more regularity conditions fail is a central issue to decide whether and to which extent to rely upon simulation.  The following simple example tries and makes the point.

\begin{es}[Testing for homogeneity in a von Mises mixture]
\label{ex:vonMises}
Suppose we observe a random sample $y_1,\ldots,y_n$ from the mixture model
\begin{equation}
\label{vMmix}
(1-p) f(y_i; 0, \kappa) + p f(y_i;\mu,\kappa), 
\end{equation}
where $0\leq p\leq 1$ is the mixing proportion.  Furthermore, $f(y_i;\mu,\kappa)$ denotes the von Mises distribution with mean direction $|\mu|\leq\pi$ and concentration parameter $\kappa\geq0$.  \cite{Fu++08} prove that the asymptotic null distribution of the likelihood ratio statistic for testing the hypothesis $p = 0$ is the squared supremum of a truncated Gaussian process.  The quantiles of the process can in principle be approximated to desirable precision by simulation, this way overcoming the algebraic difficulties of the exact solution.  However, the same authors also show that if a suitable penalisation term is used, the distribution of the corresponding modified likelihood ratio statistic converges to the simple $\chi^2_1$ distribution for $n\rightarrow\infty$.    

This is wholly different from what happens in the Gaussian case.  If the component densities $f(y_i;\mu,\kappa)$ in \eqref{vMmix} represent normal distributions with unknown mean $\mu\in{\mathbb R}$ and variance $\kappa>0$, the distribution of the likelihood ratio statistic for testing model homogeneity diverges to infinity unless suitable constraints are imposed \citep{ChenChen03}.  This is because normal mixtures with unknown variance are not identifiable unlike the von Mises mixture model~\eqref{vMmix}; see Section~\ref{gaussian-mixture}.  Trying and simulating the limiting distribution in this case would lead to totally misleading results as the likelihood ratio tends to infinity with probability one.  The finite-sample distribution of the likelihood ratio statistic, however, can be approximated using the parametric bootstrap as in \citet{McLachlan87}.  
\end{es}

The required conditions, which are typically of Cram\'er type \citep[\S 33.3]{Cramer46}, require, among others, differentiability with respect to the parameters of the underlying joint probability or density function up to a suitable order and finiteness of the Fisher information matrix.  Models which satisfy these requirements are said to be `regular' and cover a wide range of applications.  However, there are many important cases where one or more conditions break down.  A highly cited review of nonregular problems is \cite{Smith89}; see also the discussion paper by \citet{ChengTraylor95}.  Further examples can be found in \citet[\S 3.8]{BNCox94}, \citet[\S 4.6]{Davison03} and \citet[Chapter~7]{Cox06}.  A classical example, which is traditionally used to demonstrate the failure of parametric likelihood theory, is Neyman and Scott's (1948) paradox. \nocite{NeymanScott48} 

\begin{es}[Growing number of parameters]
\label{ex:growing}
Let $(X_1,Y_1), \ldots, (X_n,Y_n)$ denote $n$ independent pairs of mutually independent and normally distributed random variables such that for each $i=1,\ldots,n$, $X_i$ and $Y_i$ have mean $\mu_i$ and common variance $\sigma^2$.  The maximum likelihood estimator of $\sigma^2$ is
$$\hat\sigma_n^2=\frac{1}{2n}\sum_{i=1}^n\{(X_i-\hat\mu_i)^2+ (Y_i-\hat\mu_i)^2\},$$ 
with $\hat\mu_i=(X_i+Y_i)/2$.  Straightforward calculation shows that, for $n\rightarrow\infty$, $\hat\sigma_n^2$ converges in probability to $\sigma^2/2$ instead of the true value $\sigma^2$.  The reason is that only a finite number of observations, in fact two, is available for estimating the unknown sample means $\mu_i$.  This violates a major requirement which underlies the consistency of the maximum likelihood estimator, namely that the uncertainty of all parameter estimates goes to zero.
\end{es}

Example~\ref{ex:growing} is an early formulation of an incidental parameters problem.  Other examples of this type are reviewed in \cite{Lancaster00}, who also discusses the relevance of the Neyman--Scott paradox in statistics and economics.  A recent contribution is \cite{Feng++12}.  Non-regularity may also arise when the parameter value under the null hypothesis is not an interior point of the parameter space, or when some of the parameters disappear under the null hypothesis.  The following simple example shows what may happen when the support of the distribution depends on the parameter $\theta$. 

\begin{figure}[t]
\centering
\includegraphics[width=\linewidth]{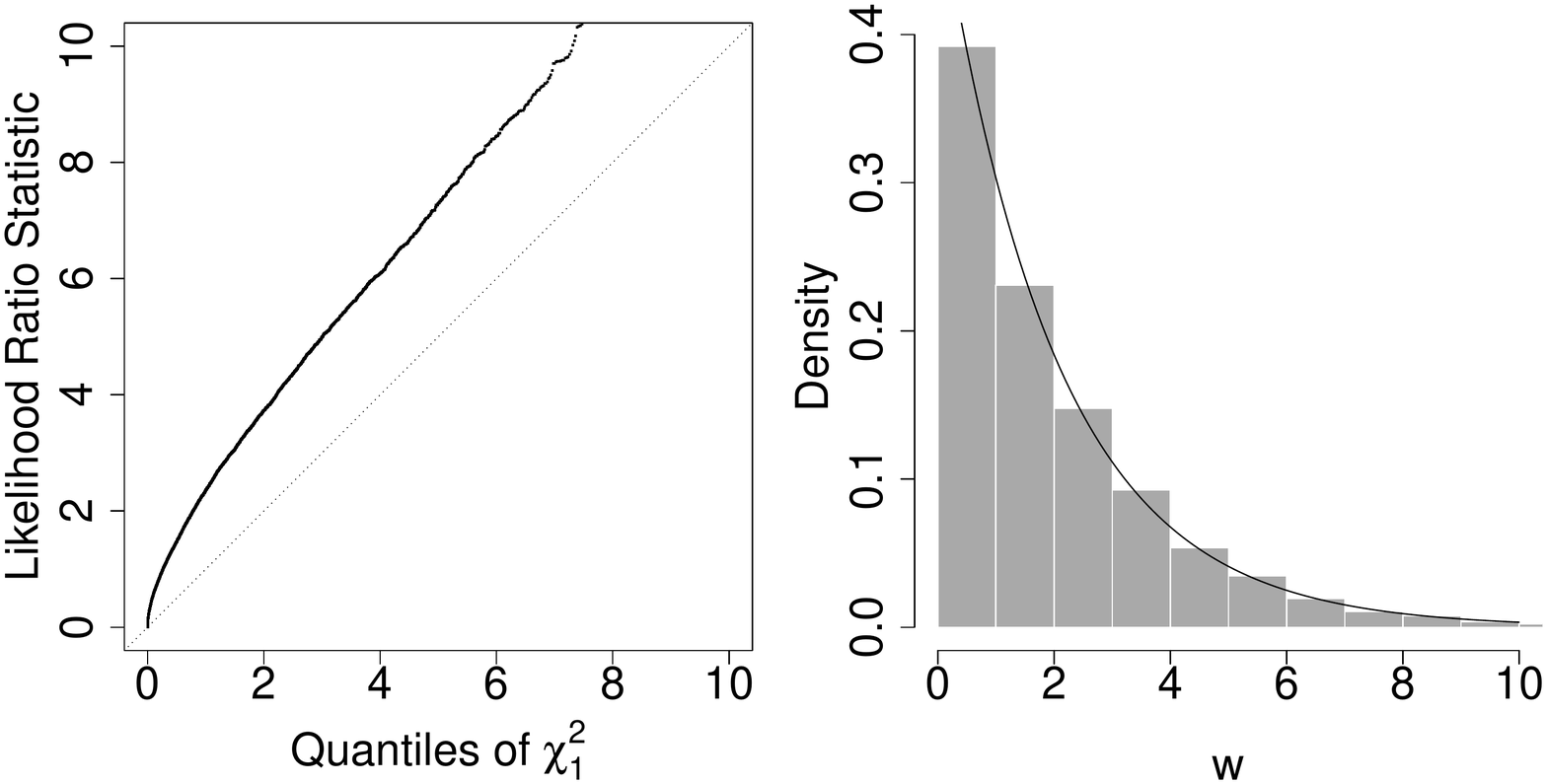}
\caption{\label{fig:ex_2} Example~\ref{ex:translated}: Translated exponential distribution.  Values of the likelihood ratio $W(3)$ observed in 10,000 exponential samples of size $n=50$ generated with rate equal to 1 and translated by $\theta_0=3$.  Left: $\chi^2_1$ quantile plot.  The diagonal dotted line is the theoretical $\chi^2_1$ approximation.  Right: histogram and superimposed $\chi^2_2$ density (solid line).}
\end{figure}

\begin{es}[Translated exponential distribution]
\label{ex:translated}
Let $X_1,\ldots,X_n$ be an independent and identically distributed sample from an exponential distribution with rate equal to 1.  Consider the translation $Y_i=X_i+\theta$, with $\theta > 0$ unknown.  Given the minimum observed value $Y_{(1)}$, the likelihood ratio statistic for testing the hypothesis that $\theta=\theta_0$ is $W(\theta_0)=2n(Y_{(1)}-\theta_0)$.  Straightforward calculation proves that under the null hypothesis $W(\theta_0)$ has a $\chi^2_2$ distribution, not the classical $\chi^2_1$ limiting distribution.  Furthermore, the maximum likelihood estimator of $\theta$ is no longer asymptotically normal.  Indeed, it is easy to show that $Y_{(1)}-\theta$ follows exactly an exponential distribution with rate $n$.  The left panel of Figure~\ref{fig:ex_2} shows the $\chi^2_1$ quantile plot of the likelihood ratio statistic observed in 10,000 exponential samples of size $n=50$ generated with rate equal to 1 and translated by $\theta_0=3$.  The finite-sample distribution of $W(3)$ is visibly far from the theoretical $\chi^2_1$ approximation represented by the dotted diagonal line.  The right panel reports the empirical distribution of the likelihood ratio statistics with superimposed the $\chi^2_2$ density (solid line).  
\end{es}

These situations are not mere mathematical artifacts, but include many models of practical interest, such as mixture distributions and change-point problems, in genetics, reliability, econometrics, and many other fields.  There is, indeed, a rich literature on this topic.  The majority of existing results consider the failure of one condition at a time, but failure of two assumptions simultaneously has also received attention.  In the absence of a unifying theory, most of the individual problems have been treated on their own.  After careful consideration, we decided to group them into three broad classes:  boundary problems, indeterminate parameter problems and change-point problems.  We furthermore restrict our attention to the key results.  Figure~\ref{timeline} depicts a personal selection of these.  The corresponding prototype derivations are provided in Appendix~A, while Appendix~B of the Supplementary Material lists further contributions, such as analogies with alternative test statistics and/or nonparametric and semiparametric models.

\begin{figure*}[t]
\centering
\label{timeline} 
\includegraphics[width=1\linewidth]{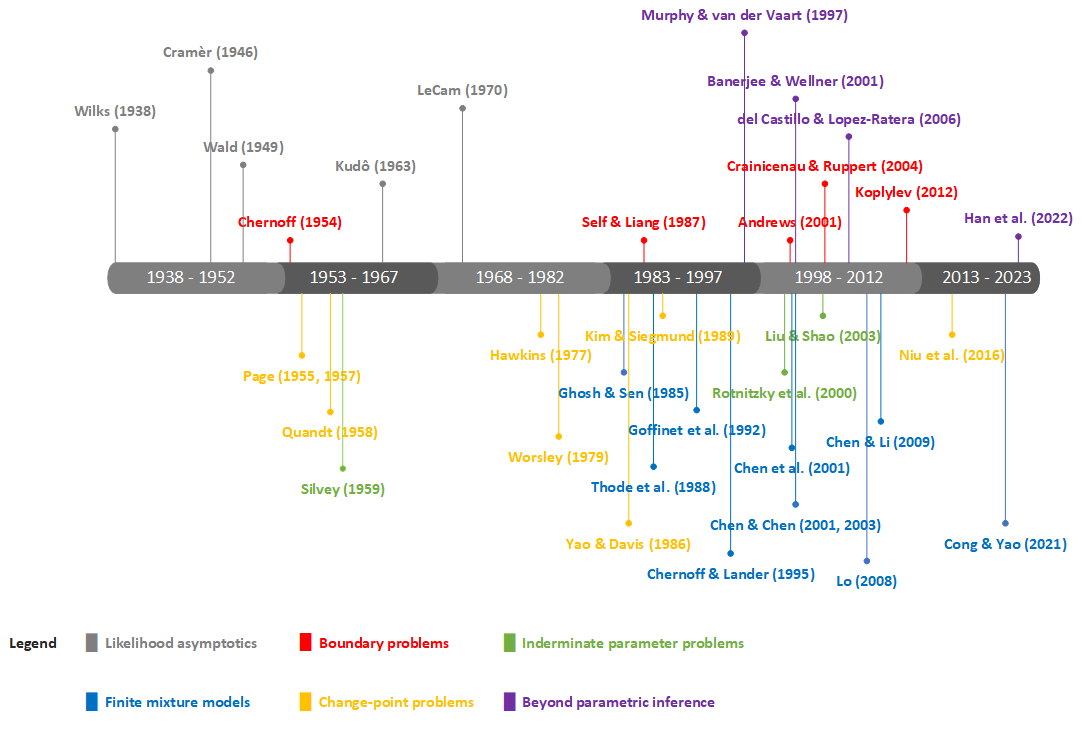}
\caption{Publishing timeline of a selection of contributions on the large-sample properties of the likelihood ratio statistic in nonregular settings.}
\end{figure*}

The paper is organised as follows.  First order parametric inference based on the likelihood function of a regular model is reviewed in Section~\ref{likelihood-asymptotics}, together with the conditions upon which it is based.  Section~\ref{boundary-problems} treats the first nonregular setting and embraces, in particular, testing for a value of the parameter which lies on its boundary.  Section~\ref{indeterminate-parameter-problems} concerns models where one part of the parameter vanishes when the remaining one is set to a particular value.  The best-studied indeterminate parameter problems are finite mixture models.  Given their widespread use in statistical practice, and their closeness to boundary problems, we will consider them separately in Section~\ref{finite-mixture-models}.  The third broad class of nonregular models, that is, change-point problems, are reviewed in Section~\ref{change-point-problems}.  Most articles investigate the consequences of the failure of one regularity condition at a time.  Mixture distributions and change-point problems deserve special attention as they represent situations where two conditions fail simultaneously.  Section~\ref{beyond-parametric-inference} reviews cases which do not fit into the above three broad model classes, but still fall under the big umbrella of nonregular problems.  These include, among others, shape constrained inference, a genre of nonparametric problem which leads to highly nonregular models.  

Despite the many remarkable theoretical developments in likelihood-based asymptotic theory for nonregular parametric models, one may wonder why the corresponding results are little known especially among practitioners.  We believe there are at least two reasons.  The first is that the results are highly scattered, in time and scope, which makes it difficult to get the general picture.  The second reason is that the limiting distributions are often fairly complex in their derivation and implementation.  Section~\ref{software} reviews the few software implementation we are aware of.  The paper then closes with the short summary discussion of Section~\ref{discussion}.


\section{Likelihood Asymptotics}
\label{likelihood-asymptotics}

\subsection{First order theory}
\label{first-order-theory}

\subsubsection{General notation.~}
\label{general-notation}
Consider a parametric statistical model $\mathcal{F} = \{f(y;\theta), \Theta, \mathcal{Y}\}$, where $y=(y_1,\ldots,y_n)$ are $n$ observations from the $d$-dimensional random variable $Y=(Y_1,\ldots,Y_n)$, $d\geq 1$, with probability density or mass function $f(y;\theta)$ whose support is $\mathcal{Y} \subseteq \mathbb{R}^d$.  Furthermore, the $p$-dimensional parameter $\theta$ takes values in a subset $\Theta\subseteq\mathbb{R}^p$, $p\geq1$.  Throughout the paper we will consider $y$ an independent and identically distributed random sample unless stated differently.  Furthermore, there may be situations where the model is specified by a function different from $f(y;\theta)$ such as the cumulative distribution function $F(y;\theta)$.  Let $L(\theta)=L(\theta;y)\propto f(y;\theta)$ and $l(\theta)=\log{L(\theta)}$ denote the likelihood and the log-likelihood functions, respectively.  The large interest in this inferential tool is motivated by the idea that $L(\theta)$ will be larger for values of $\theta$ near the true value $\theta^0$ of the distribution which generated the data.  The maximum likelihood estimate (MLE) $\hat\theta$ of $\theta$ is the value of $\theta$ which maximises $L(\theta)$ or equivalently $l(\theta)$.  That is, it is the value to which the data assign maximum evidence.  Under mild regularity conditions on the log-likelihood function to be discussed in Section~\ref{regularity-conditions}, $\hat\theta$ solves the score equation $u(\theta)=0$, where $u(\theta)=\partial l(\theta)/\partial\theta$ is the score function.  We furthermore define the observed information function $j(\theta)=-\partial^2 l(\theta)/\partial\theta\partial\theta^\top$, where $\theta^\top$ denotes transposition of $\theta$, and the expected or Fisher information $i(\theta)=E\left[j(\theta;Y)\right]$.

\subsubsection{No nuisance parameter.~}
\label{scalar_case}
The three classical likelihood-based statistics for testing $\theta=\theta_0$ are the \\[1ex]
\indent\indent
standardized MLE, \quad 
$(\hat{\theta}-\theta_0)^\top j(\hat\theta)(\hat{\theta}-\theta_0)$, \\[1ex]
\indent\indent
score statistic, \quad \quad \quad
$~u(\theta_0)^\top j(\hat\theta)^{-1}u(\theta_0)$, \\[1ex]
\indent\indent
likelihood ratio  \quad \quad \quad
$W(\theta_0)=2\{l(\hat{\theta})-l(\theta_0)\}$, \\[1.5ex]
where the observed information $j(\hat\theta)$ is at times replaced by the Fisher information $i(\theta)$.  These statistics are also known under the names of Wald's, Rao's and Wilks' tests, respectively.  If the parametric model is regular, the finite-sample null distribution of the above three statistics converges to a $\chi^2_p$ distribution to the order $O(n^{-1})$ as $n\rightarrow\infty$.  For $\theta$ scalar, inference may be based on the corresponding signed versions, that is, on the signed 
Wald statistic, 
$(\hat{\theta}-\theta_0)j(\hat\theta)^{1/2}$, 
score statistic, 
$u(\theta_0)j(\theta_0)^{-1/2}$, and 
likelihood root, 
$$r(\theta_0) = {\rm sign}(\hat{\theta}-\theta_0)[2\lbrace l(\hat{\theta})-l(\theta_0)\rbrace]^{1/2},$$ 
whose finite-sample distributions converge to the standard normal distribution to the order $O(n^{-1/2})$.

\subsubsection{Nuisance parameters.~}
Suppose now that the parameter $\theta=(\psi,\lambda) \in \Psi \times \Lambda$ is partitioned into a $p_0$-dimensional parameter of interest, $\psi\in\Psi\subseteq\mathbb{R}^{p_0}$, and a vector of nuisance parameters $\lambda\in\Lambda\subseteq\mathbb{R}^{p-p_0}$ of dimension $p-p_0$.  Large-sample inference for $\psi$ is commonly based on the profile log-likelihood function 
$$l_{\rm p}(\psi)=\sup_{\lambda\in\Lambda}l(\psi,\lambda),$$
which maximises the log-likelihood $l(\psi,\lambda)$ with respect to $\lambda$ for fixed $\psi$.  The profile likelihood ratio statistic for testing $\psi \in \Psi_0$ is
\[
W_{\rm p} (\psi_0) = 2\{\sup_{\psi\in\Psi}l_{\rm p}(\psi) - \sup_{\psi\in\Psi_0}l_{\rm p}(\psi)\},
\]
where $\Psi_0\subset\Psi$ is the parameter space specified under the null hypothesis.  If the null hypothesis is $\psi=\psi_0$, the finite-sample distribution of $W_{\rm p} (\psi_0)$  converges to the $\chi^2_{p_0}$ distribution to the order $O(n^{-1})$ for $n\rightarrow\infty$.  

If there exists a closed form expression for the constrained maximum likelihood estimate $\hat{\lambda}_{\psi}$ of $\lambda$ for given $\psi$, the profile log-likelihood function may be written as
\begin{equation}
\label{closed-lp}
l_{\rm p}(\psi)=\sup_{\lambda\in\Lambda}l(\psi,\lambda)=l(\psi,\hat{\lambda}_{\psi}).
\end{equation}
A typical situation where $\hat\lambda_\psi$ is not available in closed form is when the nuisance parameter $\lambda$ vanishes under the null hypothesis, as will be addressed in Section~\ref{non-identifiable-parameters}.  If \eqref{closed-lp} holds, we may define the profile Wald, score and likelihood ratio statistics for testing $\psi=\psi_0$ as in Section~\ref{scalar_case}, but now in terms of the profile log-likelihood $l_{\rm p}(\psi)$, with $u_{\rm p}(\psi)=\partial l_{\rm p}(\psi)/\partial\psi$ and $j_{\rm p}(\psi)=\partial l_{\rm p}(\psi)/\partial\psi\partial\psi^\top$ being the profile score and profile observed information functions.  The asymptotic null distribution of these statistics is a $\chi^2_{p_0}$ distribution up to the order $O(n^{-1})$.  If $\psi$ is scalar, the distributions of the corresponding signed versions, $(\hat{\psi}-\psi_0)j_{\rm p}(\hat\psi)^{1/2}$, $u_{\rm p}(\psi_0)j_{\rm p}(\psi_0)^{-1/2}$, and 
\begin{equation}
\label{rp}
r_{\rm p}(\psi_0)={\rm sign}(\hat{\psi}-\psi_0)[2\lbrace l_{\rm p}(\hat{\psi})-l_{\rm p}(\psi_0)\rbrace]^{1/2},
\end{equation}
may be approximated by standard normal distributions up to the order $O(n^{-1/2})$.

\subsection{Regularity conditions}
\label{regularity-conditions}

\subsubsection{Definition}
The first step in the derivation of the large-sample approximations and statistics of Sections~\ref{first-order-theory} is typically Taylor series expansion of the log-likelihood function $l(\theta)$, or quantities derived thereof, in $\hat\theta$ around $\theta$.  We illustrate this by considering the derivation of the asymptotic distribution of the likelihood ratio statistic $W(\theta)=2\{l(\hat\theta)-l(\theta)\}$ for the scalar parameter case. 

\begin{es}[Asymptotic distribution of the likelihood ratio]
\label{ex:asymptotic}
Let $p=1$ and $l_m=l_m(\theta)=d^ml(\theta)/d\theta^m$ be the derivative of order $m=2,3,\ldots$ of $l(\theta)$, the log-likelihood function for $\theta$ in a regular parametric model.  Recall that $-l_2(\theta;y)=j(\theta)$ represents the observed information, while $E[-l_2(\theta;Y)]=i(\theta)$ is the expected Fisher information.  Taylor series expansion of $l(\theta)$ around $\hat\theta$ yields
\begin{eqnarray}
\nonumber
l(\theta) & = & l(\hat\theta) - \frac{1}{2}j(\hat\theta)(\hat\theta-\theta)^2 - \frac{1}{6}(\hat\theta-\theta)^3l_3(\tilde\theta),
\end{eqnarray}
where $\tilde\theta$ is such that $|\tilde\theta-\hat\theta| < |\theta-\hat\theta|$.  Suitable rearrangement of the terms, leads to 
\begin{eqnarray}
\nonumber
W(\theta) & = & j(\hat\theta)(\hat\theta-\theta)^2 + \frac{1}{3}l_3(\tilde\theta)(\hat\theta-\theta)^3 \\
\label{expansion.lr}
          & = & \frac{j(\hat\theta)}{i(\theta)}i(\theta)(\hat\theta-\theta)^2 + \frac{1}{3}l_3(\tilde\theta)(\hat\theta-\theta)^3 . 
\end{eqnarray}
Now, under suitable regularity conditions on $l(\theta)$ and of its first three derivatives, $\hat\theta \overset{p}{\rightarrow} \theta$ \citep[Statement (i) of Theorem on p.~145]{Serfling80},
\begin{eqnarray}
\nonumber
\frac{j(\theta)}{i(\theta)} \overset{p}{\rightarrow} 1 \quad \text{and} \quad
i(\theta)^{1/2}(\hat\theta-\theta) \overset{d}{\rightarrow} Z, 
\end{eqnarray}
where $Z$ has the standard normal distribution \citep[Lemma~B and Lemma~A(ii)]{Serfling80}.  Furthermore, by the law of large numbers $l_3(\theta) \overset{p}{\rightarrow} c < +\infty$.  Advocating Slutzky’s lemma, the leading term in (\ref{expansion.lr}) hence converges asymptotically to the $\chi^2_1$ distribution, while the second addend is of order $o_p(1)$.  This leads to the well known result for Wilks' statistic.  
\end{es}

The derivation of Example~\ref{ex:asymptotic} requires that the model under consideration is regular.  This implies first of all that the log-likelihood function can be differentiated at least to the third order, but also that the expected values of log-likelihood derivatives are finite and that their asymptotic order is proportional to the sample size.  \cite{Wald49}---who is generally acknowledged for having provided the earliest proof of consistency of the maximum likelihood estimator which is mathematically correct---furthermore emphasized the importance of the compactness of the parameter space $\Theta$ and that the maximum likelihood estimator be unique.  Indeed, the former condition was missing in Cram\'er's (1946) and Huzurbazar's (1948) \nocite{Huzurbazar48} proofs.      

In this paper, by the term ``regularity conditions’’ we mean the assumptions on the parametric statistical model $\mathcal{F}$ that ensure the validity of classical asymptotic theory.  These may be formulated in several ways; see e.g.\ \citet[p.~281]{CoxHinkley74}, \citet[\S 3.8]{BNCox94}, \citet[\S 3.2.3]{Azzalini96}, \citet[\S 4.7]{Severini00}, \citet[Chap.~5]{vanderVaart00}, \citet[\S 4.6]{Davison03}, \citet[\S 6.1, \S 6.2 and A.1]{Hogg++19}.  We will assume that the following five conditions on $\mathcal{F}=\{f(y;\theta), \Theta, \mathcal{Y}\}$ and related likelihood quantities hold.
\begin{description}
\item[{\it Condition~1}] All components of $\theta$ are identifiable. That is, two probability density or mass functions $f(y;\theta^1)$ and $f(y;\theta^2)$ defined by any two different values $\theta^1\neq\theta^2$ of $\theta$ are distinct almost surely.  
\end{description}

\begin{description}
\item[{\it Condition~2}] The support $\mathcal{Y}$ of $f(y;\theta)$ does not depend on $\theta$.  
\end{description}

\begin{description}
\item[{\it Condition~3}] The parameter space $\Theta$ is a compact subset of $\mathbb{R}^p$, for a fixed positive integer $p$, and the true value $\theta^0$ of $\theta$ is an interior point of $\Theta$.  
\end{description}

\begin{description}
\item[{\it Condition~4}] The partial derivatives of the log-likelihood function $l(\theta;y)$ with respect to $\theta$ up to the order three exist in a neighbourhood of the true parameter value $\theta^0$ almost surely.  Furthermore, in such a neighbourhood, $n^{-1}$ times the absolute value of the log-likelihood derivatives of order three are bounded above by a function of $Y$ whose expectation is finite.  
\end{description}

\begin{description}
\item[{\it Condition~5}] The first two Bartlett identities hold, which imply that 
\[E[u(\theta;Y)] = 0 , \quad i(\theta) = {\rm Var}[u(\theta;Y)],\]
in addition to $ 0 < {\rm Var}[u(\theta;Y)] < \infty$.
\end{description}

\subsubsection{Failure of regularity}
Conditions~1--5 are relevant in many important models of practical interest, and can fail in as many ways.  For instance, from the perspective of significance testing, Condition~1 fails when under the null hypothesis parameters defined for the whole model become undefined and therefore inestimable.  We already mentioned this situation when introducing the profile log-likelihood function; non-identifiability of the parameters will be addressed in Section~\ref{non-identifiable-parameters}.  Further examples are treated in Sections~\ref{singular-information-matrix} and \ref{finite-mixture-models}.  Failure of Condition~2 is addressed in \cite{HiranoPorter03} and \cite{Severini04}.  Failure of Condition~3 characterises the first and most extensively explored nonregular setting, that is, boundary problems; see Section~\ref{boundary-problems}.  Furthermore, they include the, to our knowledge, only contribution which explores the higher order properties of likelihood-based test statistics in a nonregular setting \citep{CastilloLopez06}.  The compactness condition, in particular, can be omitted, provided it is replaced by some other requirements; see, for instance, \citet[Page~119]{Pfanzagl17}.  This will be also the case for a number of the large-sample results derived for nonregular models; see, for instance, Section~\ref{finite-mixture-models}.  Condition~4 typically does not hold in change-point problems, which will be treated in Section~\ref{change-point-problems}.  A further prominent example where Condition~4 is not satisfied, is the double exponential, or Laplace, distribution, which arises in quantile regression.  For a book-length review of this topic we refer the reader to \cite{Koenker++17}.  Condition~5 is guaranteed if standard results on the interchanging of integration and differentiation hold, Condition~2 is satisfied, and the log-likelihood derivatives are continuous functions of $\theta$.  A typical situation where this condition fails is when the data under analysis are derived from a probability density which does not belong to the model $f(y;\theta)$, a topic of much investigation in robustness \citep{HuberRonchetti09}.  A remedy is provided by Godambe's theory of estimating equations \citep{Godambe91}. 

Conditions~4 and 5, as used by \cite{Cramer46}, \cite{Wald49} and others, require the existence of at least three derivatives of the log-likelihood function together with some uniform integrability restrictions.  Condition~4, in particular, embraces both, the existence of the partial derivatives of $l(\theta)$ and their asymptotic order.  An example for which this latter condition does not hold is the Pearson Type III (or translated Gamma) distribution \citep{Blischke++69}, which generalizes Example~\ref{ex:translated}.  In this latter case, $|d l(\theta;y)/d\theta|=n$ is not dominated by an integrable function on $(\theta, +\infty)$.  Condition~4 ensures the consistency of the maximum likelihood estimator $\hat\theta$, and the existence of a quadratic approximation to the likelihood ratio function $l(\hat\theta) - l(\theta)$ in the Euclidean $n^{-1/2}$-neighbourhood of $\theta$ of the form
\begin{eqnarray}
\nonumber
l(\hat\theta) - l(\theta) & = &  (\hat\theta-\theta)^\top u(\theta;y) \\
\label{quadratic}
             & - & \frac{1}{2}(\hat\theta -\theta)^\top i(\theta)(\hat\theta-\theta) + o_p(1), 
\end{eqnarray}
which involves the score function $u(\theta;y)$ and the expected Fisher information $i(\theta)$.  However, these conditions do not have by themselves any direct statistical interpretation.  \cite{LeCam70} presents a different type of regularity assumption---differentiability in quadratic mean---which requires only differentiation of order one and may be justified from a statistical point of view.

\subsection{Local asymptotics}
\label{local-asymptotics}
To give a glimpse of LeCam’s ideas, assume, without loss of generality, $\theta\in\mathbb{R}$ scalar.  A statistical model $\mathcal{F} = \{f(y_i;\theta), \theta\in\Theta\}$ is said to be {\it differentiable in quadratic mean} (DQM) at $\theta$ if there exists a random function $\dot{u}(\theta;Y_i)$ such that, as $h \rightarrow 0$,
\begin{eqnarray}
\nonumber
E_{\theta} \left[ \left\{ \sqrt{f(Y_i;\theta+h)} -\sqrt{f(Y_i;\theta)} \right . \right .\\
\label{DQMcondition}
\left. \left. -\frac{h}{2} \dot{u}(\theta;Y_i)\sqrt{f(Y_i;\theta)} \right\}^2 \right] = o(||h^2||).
\end{eqnarray} 
The expectation is taken with respect to the distribution indexed by $\theta$ and the function $\dot{u}(\theta;Y_i)$ is said to be the quadratic mean derivative of the square root $\sqrt{f(Y_i;\theta)}$ of the probability density function \citep[Chapter~12]{LehmannRomano05}.  As shown in LeCam’s (1970) paper, the regularity conditions of Cram\'er type imply differentiability in quadratic mean.  Indeed, for every differentiable $f(y_i;\theta)$ 
\begin{equation}
\label{DQM}
\frac{\partial}{\partial\theta} \sqrt{f(y_i;\theta)} = \frac{1}{2} u(\theta;y_i) \sqrt{f(y_i;\theta)},
\end{equation}
where $\dot{u}(\theta;y_i)=u(\theta;y_i)$ is the score function as defined in Section~\ref{first-order-theory}.  The opposite does not hold true, a prominent counterexample being the Laplace distribution.  Differentiability in quadratic mean hence generalizes Condition~5 in a natural way as it is possible to show that $E_{\theta}[\dot{u}(\theta;Y_i)] = 0$ and that the equivalent of the unit expected Fisher information $\ddot{\imath}_1(\theta) = E_{\theta}[\{\dot{u}(\theta;Y_i)\}^2]$ is finite.  

Using differentiability in quadratic mean, LeCam gives rise to a radically different type of asymptotic inference called \textit{local asymptotics}.  The word `local’ is meant to indicate that one looks at a sequence of alternative hypotheses of the form $\theta_n = \theta + \epsilon/\sqrt{n}$, where $\epsilon$ is any given real number.  The properties of the likelihood-based procedures are hence studied in a small neighbourhood $\theta \pm \epsilon/\sqrt{n}$ of the fixed parameter $\theta$ defined by $\epsilon$, where `small’ means of size $O(1/\sqrt{n})$.  The motivation for studying a local approximation is that, usually, asymptotically, the `true’ parameter value can be known with unlimited precision. The real difficulty is therefore to distinguish between values which are `close’.  Closeness in this case is measured in terms of the Hellinger distance
\[
H^2(\theta) = \frac{1}{2}E_\theta\left[ \left\{\sqrt{\lambda(\theta;Y_i)}-1\right\}^2 \right],
\]
with 
\[
\lambda(\theta;Y_i) = \frac{f(Y_i;\theta+h)}{f(Y_i;\theta)},
\]
whose definition can be linked to the notion of differentiability in quadratic mean as
\begin{eqnarray*}
H^2(\theta) & = & \frac{h^2}{8} E_\theta \left[ \left\{\dot{u}(\theta;Y_i) \right\}^2 \right] + o(||h||^2) \\
            & = & \frac{h^2}{8} \ddot{\imath}(\theta) + o(||h||^2)
\end{eqnarray*}
if the model satisfies the DQM condition \eqref{DQMcondition}.  Indeed, in this case it is possible to show that the likelihood ratio function of a random sample, $y=(y_1,\ldots,y_n)$, of size $n$ 
\[ 
lr(\theta;y) = \sum_{i=1}^n \log \lambda(\theta;y_i) 
\]
with $h=\epsilon/\sqrt{n}$, is {\it locally asymptotically quadratic} in that 
\begin{eqnarray*}
lr(\theta;y) =
\frac{\epsilon}{\sqrt{n}} \dot{u}(\theta;y) - \frac{1}{2}\frac{\epsilon^2}{n}\ddot{\imath}(\theta) + o_p(1).
\end{eqnarray*}
Note how this expression mimics the quadratic approximation \eqref{quadratic} of classical likelihood-based asymptotics, where $\hat\theta_n - \theta = O_p(1/\sqrt{n})$ and the score and the expected Fisher information functions are replaced by $\dot{u}(\theta;Y)$ and $\ddot{\imath}(\theta)=n\ddot{\imath}_1(\theta)$.  For large $n$ we can locally approximate the likelihood ratio function by the normal distribution
\[
N\left(-\frac{1}{2}\epsilon^2\ddot{\imath}_1(\theta), \ \epsilon^2\ddot{\imath}_1(\theta)\right),
\]
which then serves as the basis for the derivation of the limiting distributions of estimators and test statistics.  For further details see the two monographs by \citet[Chapter~7]{vanderVaart00} and \cite{LeCamYang90}.

In the remainder of the paper, we review the most common situations where one or some of Conditions~1--5 fail.  We will also provide some summary insight into the main prototype derivations of the corresponding finite-sample or asymptotic results.  The vast majority of the proofs require conditions of Cram\'er type; in some occasions, as for instance in Section~\ref{non-identifiable-parameters}, LeCam's local asymptotic theory will be used.


\section{Boundary Problems}
\label{boundary-problems}

\subsection{Definition}
\label{boundary-problems-definition}
\begin{figure}[tb]
\centering
\includegraphics[width=\linewidth]{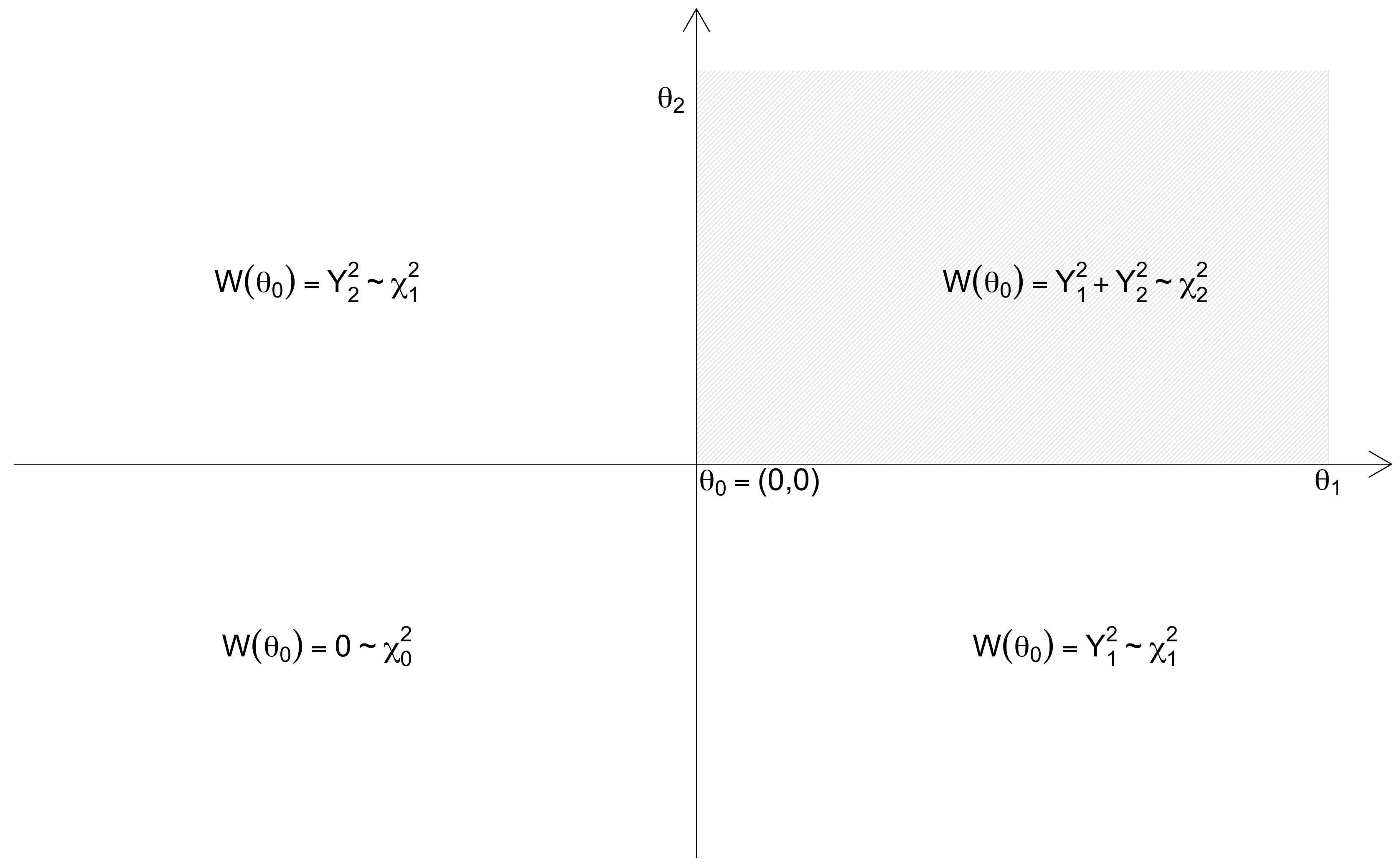}
\caption{\label{fig:1} Example~\ref{ex:boundary}: Bivariate normal.  The grey shaded area represents the parameter space $\Theta$.  Under the null hypothesis $\theta_0=(0,0)$, the parameter space collapses with the origin.  The asymptotic distribution of the corresponding likelihood ratio statistics is a mixture of $\chi^2_0$, $\chi^2_1$ and $\chi^2_2$ distributions with weights $(0.25, 0.5, 0.25)$.}
\end{figure}

A boundary problem arises when the value $\theta_0$ specified by the null hypothesis, or parts of it, is not an interior point of the parameter space.  In general terms, the ``boundary'' of the parameter space $\Theta$ is the set of values $\theta$ such that every neighbourhood of $\theta$  contains at least one interior point of $\Theta$ and at least one point which is not in $\Theta$.  Informally, the methodological difficulties in likelihood-based inference occur because the maximum likelihood estimate can only fall `on the side' of $\theta_0$ that belongs to the parameter space $\Theta$.  This implies that if the maximum occurs on the boundary, the score function need not be zero and the distributions of the related likelihood statistics will not converge to the typical normal or chi-squared distributions.  Because of the difficulties inherent the derivation of the limiting distribution of the likelihood ratio statistic, especially practitioners tend to ignore the boundary problem and to proceed as if all parameters where interior points of $\Theta$.  This is commonly called the {\it na\"\i ve} approach.  An alternative approach is to suitably enlarge the parameter space so as to guarantee that the likelihood ratio maintains the common limiting distribution; see, for instance, \cite{FengMcCulloch92}.  However, this idea works only as long as the null hypothesis is uniquely identified.  The following example gives a flavour of the statistical issues. 
\begin{es}[Bivariate normal]\label{ex:boundary}
Consider a single observation $y=(y_1,y_2)$ from the bivariate normal random variable $Y = (Y_1, Y_2)\sim N_2(\theta, I_2)$, where $\theta=(\theta_1,\theta_2)$, with $\theta_1\geq0$ and $\theta_2\geq0$, and $I_2$ is the $2\times2$ identity matrix.  Straightforward calculation shows that the null distribution of the likelihood ratio statistic for $\theta_0=(0,0)$ versus the alternative hypothesis that at least one equality does not hold, converges to a mixture of a point mass $\chi^2_0$ at 0 and two chi-squared distributions, $\chi^2_1$ and $\chi^2_2$ \cite[Example~21.3]{DasGupta08}.  Figure~\ref{fig:1} provides a graphical representation of the problem.  Because of the limitedness of the parameter space, we have that $\hat\theta_1 = \max(y_1,0)$ and $\hat\theta_2=\max(y_2,0)$.  The grey shaded area is the parameter space into which the MLE is bound to fall.  However, the random observation $Y = (Y_1, Y_2)$ can fall into any of the 4 quadrants of ${\mathbb R}^2$ with equal probability $1/4$.  When $Y$ falls into the first quadrant, that is, when $y_1, y_2 >0$, the likelihood ratio statistic is $W(\theta_0)= Y_1^2+Y_2^2$ and follows the common $\chi^2_2$ distribution.  However, if $y_1>0$ and $y_2<0$ or when $y_1<0$ and $y_2>0$, we have that $W(\theta_0)=Y_1^2\sim \chi^2_1$ and $W(\theta_0)=Y_2^2\sim\chi^2_1$, respectively.  Lastly, when $Y$ lies in the third quadrant, $W(\theta_0)= 0$ and its distribution is a point mass in 0.  Summing up, we can informally write
\begin{equation}
\label{example4}
W(\theta_0) \sim \frac{1}{4}\chi^2_0 + \frac{1}{2}\chi^2_1 + \frac{1}{4}\chi^2_2.
\end{equation}
\end{es} 

Distribution (\ref{example4}) is a special case of the so-called chi-bar squared distribution \citep{Kudo63}, denoted by $\bar\chi^2(\omega,N)$, with cumulative distribution function  
\[
\label{chi-bar} 
{\rm Pr}(\bar\chi^2 \leq c) = \sum_{\nu=0}^{N}\omega_\nu{\rm Pr}(\chi^2_\nu \leq c), 
\]
which corresponds to a mixture of chi-squared distributions with degrees of freedom $\nu$ from 0 to $N$.  In some cases, explicit and computationally feasible formulae are available for the weights $\omega=(\omega_0,\ldots,\omega_N)$.  Extensive discussion on their computation and use, with special emphasis on inequality constrained testing, is given in \citet[Chapters 2 and 3]{Robertson++88}, \cite{Wolak87}, \cite{Shapiro85, Shapiro88} and \cite{Sun88}.

\subsection{General results}\label{boundary_general_results}
The research on boundary problems was initiated by \cite{Chernoff54} who derives the asymptotic null distribution of the likelihood ratio statistic for testing whether $\theta$ lies on one or the other side of a smooth $(p-1)$-dimensional surface in a $p$-dimensional space when the true parameter value lies on the surface.  Using a geometrical argument, Chernoff established that this distribution is equivalent to the distribution of the likelihood ratio statistic for testing suitable restrictions on the mean of a multivariate normal distribution with covariance matrix given by the inverse of the Fisher information matrix using a single observation.  In particular, Chernoff proved that the limiting distribution is a $\bar\chi^2(\omega,1)$ distribution, with $\omega=(0.5,0.5)$, that is, a mixture of a point mass at zero and a $\chi^2_1$, with equal weights.  This generalizes \cite{Wilks38} result when the parameter space under the null hypothesis is not a hyperplane. 

The no doubt cornerstone contribution which inspired many researchers and fuelled an enormous literature, is the highly-cited article by \cite{SelfLiang87}.  In \cite{Chernoff54}, the parameter spaces $\Theta_0$ and $\Theta_1$, specified by the null and the alternative hypotheses, are assumed to have the same dimension.  Furthermore, the true parameter value falls on the boundary of both, $\Theta_0$ and $\Theta_1$, while it is still an interior point of the global parameter space $\Theta=\Theta_0\cup\Theta_1$.  Using geometrical arguments similar to those of \cite{Chernoff54}, \cite{SelfLiang87} study the asymptotic null distribution of the likelihood ratio statistic for testing the null hypothesis $\theta\in \Theta_0$ against the alternative $\theta \in \Theta_1=\Theta\setminus\Theta_0$.  This time, the true parameter value $\theta^0$ no longer needs be an interior point, but can fall onto the boundary of $\Theta$.  The two sets $\Theta$ and $\Theta_0$ must be regular enough to be approximated by two cones, $C_{\Theta}$ and $C_{\Theta_0}$, with vertex at $\theta_0$ \cite[Definition~2]{Chernoff54}.  Under this scenario and provided their Assumptions~1--4 hold---which translate into our Conditions~1--2 and 4--5, with likelihood derivatives taken from the appropriate side---\citet[Theorem~3]{SelfLiang87} show that the distribution of the likelihood ratio converges to the distribution of
\begin{eqnarray}
\label{selfliang}
\sup_{\theta\in C_{\Theta-\theta^0}}\left\lbrace - (\tilde Z-\theta)^\top i_1(\theta^0)(\tilde Z-\theta)\right\rbrace - \\
\nonumber
\sup_{\theta\in C_{\Theta_0-\theta^0}}\left\lbrace - (\tilde Z-\theta)^\top i_1(\theta^0)(\tilde Z-\theta)\right\rbrace .
\end{eqnarray}
Here, $C_{\Theta-\theta^0}$ and $C_{\Theta_0-\theta^0}$ are the translations of the cones $C_{\Theta}$ and $C_{\Theta_0}$, such that their vertices are at the origin, and $\tilde Z$ is a multivariate Gaussian variable with mean 0 and covariance matrix given by $i_1(\theta^0)^{-1}$, which is the Fisher information matrix for a single observation.  If we transform the random variable $\tilde Z$ so that it follows a multivariate standard Gaussian distribution $Z$, we can re-express Equation~\eqref{selfliang} as
\begin{eqnarray}
\label{selfliangstar}
\inf_{\theta\in \tilde{C}_{0}} ||Z-\theta||^2 -\inf_{\theta\in \tilde{C}} ||Z-\theta||^2 = \nonumber \\
||Z-\mathcal{P}_{\tilde{C}_0}(Z)||^2 - ||Z-\mathcal{P}_{\tilde{C}}(Z)||^2,
\end{eqnarray}
where $\tilde{C}$ and $\tilde{C}_0$ are the corresponding transformations of the cones $C_{\Theta-\theta^0}$ and $C_{\Theta_0-\theta^0}$ and $||\cdot||$ is the Euclidean norm.  Finding the null distribution requires to work out the two projections $\mathcal{P}_{\tilde{C}}(Z)$ and $\mathcal{P}_{\tilde{C}_0}(Z)$ of $Z$ onto the cones $\tilde{C}$ and $\tilde{C}_0$.  This must be done on a case by case basis as shown by the following revisitation of Example~\ref{ex:boundary}.

\begin{es}[Bivariate normal revisited]
\label{ex:boundary2}
In Example~\ref{ex:boundary} we faced a typical non-standard situation where both components of the parameter $\theta$ are of interest and both lie on the boundary of the parameter space.  Here, the Fisher information matrix is the identity matrix which is why $\tilde Z = Z=Y$ and the original two set $\Theta$ and $\Theta_0$ agree with the approximating cones.  That is, the grey shaded region $[0,\infty)\times[0,\infty)$ in Figure~\ref{fig:1} represents the sets $\Theta=C_{\Theta}=C_{\Theta-\theta_0}=\tilde{C}$, while the origin $\{0\}$ corresponds to the sets $\Theta_0=C_{\Theta_0}=C_{\Theta_0- \theta_0}=\tilde{C}_0$.  The derivation of the second term of \eqref{selfliangstar} depends on the projection of $Z$ onto $\tilde{C}$, which is 
\begin{equation*}
\mathcal{P}_{\tilde{C}}(Z)=\begin{cases}
Z=(Z_1,Z_2) & {\rm if}\quad Z_1, Z_2 > 0 \\
Z_2 & {\rm if}\quad Z_1<0, Z_2>0 \\
0 & {\rm if}\quad Z_1, Z_2 < 0 \\
Z_1 & {\rm if}\quad Z_1>0, Z_2<0, \\
\end{cases}
\end{equation*}
while $\mathcal{P}_{\tilde{C_0}}(Z)=0$.  As shown in Example~\ref{ex:boundary}, $\mathcal{P}_{\tilde{C}}(Z)$ takes on the four possible values with equal probability $1/4$.  By simple algebra, we can prove that the distribution of the likelihood ratio statistics is given by the mixture of Equation~\eqref{example4}. 
\end{es} 

\cite{SelfLiang87} present a number of special cases in which the representations \eqref{selfliang} and \eqref{selfliangstar} are used to derive the asymptotic null distribution of the likelihood ratio statistic.  In most cases, the limiting distribution is a chi-bar squared distribution whose weights depend, at times in a rather tricky way, on the partition of the parameter space induced by the geometry of the cones.  A sketch of the derivation of Equation~\eqref{selfliang} is given in Appendix~\ref{ex:self_liang}.  The proof consists of two steps.  We first consider a quadratic Taylor series expansion of the log-likelihood $l(\theta)$ around $\theta^0$, the true value of the parameter.  The asymptotic distribution of the likelihood ratio statistic is then derived as in \cite{Chernoff54} by approximating the sets $\Theta$ and $\Theta_0$ using the cones $C_{\Theta}$ and $C_{\Theta_0}$.   

A further major step forward in likelihood asymptotics for boundary problems was marked by \cite{KopylevSinha11} and \cite{Sinha++12}.  Now, the null distribution of the likelihood ratio statistic is derived by using algebraic arguments.  From the technical point of view, the derivation of a closed form expression for the limiting distribution of the likelihood ratio becomes the more difficult the more nuisance parameter lie on the boundary of the parameter space.  In particular, the derivation of the limiting distribution becomes awkward when there are more than four boundary points and/or the Fisher information matrix is not diagonal.  \cite{Sinha++12} furthermore show that when one or more nuisance parameters are on the boundary, following the na\"\i ve approach can result in inferences which are anticonservative.  In general, the asymptotic distribution turns out to be a chi-bar squared distribution with weights that depend on the number of parameters of interest and of nuisance parameters, and on where these lie in $\Theta$.  However, limiting distributions other than the $\bar\chi^2$ distribution are found as well; see, for instance, Theorem~2.1 of \cite{Sinha++12}. 

A concise review of the cases considered in \cite{SelfLiang87}, \cite{KopylevSinha11} and \cite{Sinha++12}, with some interesting examples and an account of the areas of interest in genetics and biology, is given by \cite{Kopylev12}.  The following two sections treat two special cases, namely testing for a zero variance component and constrained one-sided tests.  We will mention the mainstream contributions, while further related work can be found in Appendix~B of the Supplementary Material.  This includes, for instance, alternatives which avoid the calculation of the mixing weights of the $\bar\chi^2$ distribution and/or lead to the classical $\chi^2$ limiting distribution.

\subsection{Null variance components}
\label{boundary_null_variance}
In linear and generalized linear mixed models a boundary problem arises as soon as we want to assess the significance of one or more variance components.  The two reference papers are \cite{CrainiceanuRuppert04} and \cite{StramLee94}.  Both consider a linear mixed effects model and test for a zero scalar variance component.  However, \cite{StramLee94} assume that the data vector can be partitioned into a large number of independent and identically distributed sub-vectors, which needs not hold for \cite{CrainiceanuRuppert04}.  The limiting distributions are derived from the spectral decomposition of the likelihood ratio statistic.

More precisely, assume the following model holds,
$$Y=X\beta+Zb+\varepsilon,$$ 
where $Y$ is a vector of observations of dimension $n$, $X$ is a $n\times p$ fixed effects design matrix and $\beta$ is a $p$-dimensional vector of fixed effects.  In addition, $Z$ is a $n\times k$ random effects design matrix and $b$ is a \textit{k}-dimensional vector of random effects which are assumed to follow a multivariate Gaussian distribution with mean $0$ and covariance matrix $\sigma^2_{b}\Sigma$ of order $k\times k$.  The error term $\varepsilon$ is assumed to be independent of $b$ and distributed as a normal random vector with zero mean and covariance matrix $\sigma_{\varepsilon}^2I_n$, where $I_n$ is the identity matrix.  Suppose we are interested in testing 
\begin{equation*}
H_0 : \beta_{p+1-q}= \beta^0_{p+1-q},\ldots , \beta_{p} = \beta^{0}_p,\quad \quad \sigma^2_b = 0
\end{equation*}
\textrm{against}
\begin{equation*}
H_1 : \beta_{p+1-q}\neq \beta^0_{p+1-q},\ldots , \beta_{p} \neq \beta^{0}_p , \quad\textrm{or}\quad \sigma^2_b > 0
\end{equation*}
for some positive value of $q\in\{1,\ldots,p\}$.  Non-regularity arises as under the null hypothesis $\sigma^2_b = 0$ falls on the boundary of the parameter space.  Furthermore, the alternative hypothesis that $\sigma^2_b > 0$ induces dependence among the observations $Y$.  \citet[Theorem~1]{CrainiceanuRuppert04} show that the finite-sample distribution of the likelihood ratio statistic agrees with the distribution of
\begin{equation}\label{spectral}
n\left( 1+\frac{\sum_{s=1}^qu^2_s}{\sum_{s=1}^{n-p}w^2_s}\right)+\sup_{\lambda\geq 0} f_{n}(\lambda),
\end{equation}
where $u_s$ for $s=1,\ldots,k$ and $w_s$ for $s=1,\ldots, n-p$ are independent standard normal variables, $\lambda=\sigma^2_b/\sigma^2_\varepsilon$, and
\begin{equation*}
f_n(\lambda)=n\log{\left\lbrace 1+\frac{N_n(\lambda)}{D_{n}(\lambda)}\right\rbrace }-\sum_{s=1}^k\log{\left( 1+\lambda\xi_{s,n}\right) },
\end{equation*} 
where
$$N_n(\lambda)=\sum_{s=1}^k\frac{\lambda\mu_{s,n}}{1+\lambda\mu_{s,n}}w^2_s,$$
and 
$$D_n(\lambda)=\sum_{s=1}^k\frac{w^2_s}{1+\lambda\mu_{s,n}}+\sum_{s=k+1}^{n-p}w^2_s.$$ 
Here, $\mu_{s,n}$ and $\xi_{s,n}$ are the $k$ eigenvalues of the matrices $\Sigma^{\frac{1}{2}}Z^TP_0Z\Sigma^{\frac{1}{2}}$ and $\Sigma^{\frac{1}{2}}Z^TZ\Sigma^{\frac{1}{2}}$, respectively.  The matrix $P_0=I_n-X(X^TX)^{-1}X^{T}$ is the matrix which projects onto the orthogonal complement to the subspace spanned by the columns of the design matrix $X$.  Theorem~2 of \cite{CrainiceanuRuppert04} shows that the asymptotic null distribution of the likelihood ratio statistic depends on the asymptotic behaviour of the eigenvalues $\mu_{s,n}$ and $\xi_{s,n}$.  The limiting distribution, in general, differs from the chi-bar squared distribution which often holds for independent and identically distributed data.  

Formula \eqref{spectral} represents the spectral decomposition of the likelihood ratio statistic.  A similar result is also derived for the restricted likelihood ratio \citep[Formula~9]{CrainiceanuRuppert04}.  The unquestioned advantage of these two results is that they allow us to simulate the finite-sample null distribution of the two test statistics once the eigenvalues are calculated.  Furthermore, this simulation is more efficient than bootstrap resampling, as the speed of the algorithm only depends on the number of random effects $k$, and not on the number of observations $n$.  Applications of Crainiceanu and Ruppert's (2004) results include testing for level- or subject-specific effects in a balanced one-way ANOVA, testing for polynomial regression versus a general alternative described by P-splines and testing for a fixed smoothing parameter in a P-spline regression.

\subsection{Constrained one-sided tests}
\label{constrained-one-sided-tests}
Multistage dose-response models are a further example of boundary problem.  A $K$-stage model is characterised by a dose-response function of the form 
\[ \psi(d;\beta) = \psi(\beta_0 + \beta_1 d+ \beta_2 d^2 + \cdots  + \beta_K d^K), \]
where $d$ is the tested dose and $\psi(\cdot)$ is a function of interest such as, for instance, the probability of developing a disease.  The coefficients $\beta_k \geq 0$, for $k=1,\ldots,K$, are often constrained to be non-negative so that the dose-response function will be non-decreasing.  There is no limit on the number of stages $K$, though in practice this is usually specified to be no larger than the number of non-zero doses.  Testing whether $\beta_k = 0$ results in a boundary problem and requires the application of a so-called constrained one-sided test.  Apart from clinical trials, constrained one-sided tests are common in a number of other areas, where the constraints on the parameter space are often natural such as testing for over-dispersion, for the presence of clusters and for homogeneity in stratified analyses.  All these instances amount to having the parameter value lying on the boundary of the parameter space under the null hypothesis.  Despite their importance in statistical practice, few contributions are available on the asymptotic behaviour of the most commonly used test statistics, and of the likelihood ratio in particular.

A first contribution which evaluates the asymptotic properties of constrained one-sided tests is \cite{Andrews01}, who establishes the limiting distributions of the Wald, score, quasi-likelihood and rescaled quasi-likelihood ratio statistics under the null and the alternative hypotheses.  The results are used to test for no conditional heteroscedasticity in a GARCH(1,1) regression model and zero variances in random coefficient models.  \cite{SenSilvapulle02} review refinements of likelihood-based inferential procedures for a number of parametric, semiparametric, and nonparametric models when the parameters are subject to inequality constraints.  Special emphasis is placed on their applicability, validity, computational flexibility and efficiency.  Again, the chi-bar squared distribution plays a central role in characterising the limiting null distribution of the test statistics, while the corresponding proof requires tools of convex analysis, such projections onto cones.  See \cite{SilvapulleSen05} for a book-length account of constrained statistical inference.


\section{Indeterminate parameter problems}
\label{indeterminate-parameter-problems}

\subsection{Definition}
An ``indeterminate parameter’’ problem occurs when setting one of the components of the parameter $\theta=(\theta_1,\theta_2)$ to a particular value, say $\theta_1 = \theta_{10}$, leads to the disappearance of some or all components of $\theta_2$.  The model is no longer identifiable, as all probability density or mass functions $f(y;\theta)$ with $\theta_1=\theta_{10}$ and arbitrary $\theta_2$ identify the same distribution.  The following simple example illustrates this point. 
\begin{es}[Loss of identifiability in jump regression] 
\label{jump-regression}
Consider the model
\[
Y=\theta_{11}+\theta_{12} \mathbbm{1}(X>\theta_2)+\varepsilon,\quad \varepsilon\sim f(\varepsilon),
\] 
where $Y$ is a continuous response, $X$ a corresponding covariate and $\mathbbm{1}(X>\theta_2)$ represents the indicator function which assumes value 1 if $X>\theta_2$ and zero otherwise.  Furthermore, $\theta_1=(\theta_{11},\theta_{12})$ is a real valued vector of regression coefficients, while $\theta_2\in\mathbb{R}$ defines the point at which the jump occurs.  Assume that $\varepsilon$ is a zero-mean error term with density function $f(\varepsilon)$.  The mean of the variable $Y$ is $\theta_{11}$ for values of $X$ less or equal to $\theta_2$ and is equal to $\theta_{11}+\theta_{12}$ for values of $X$ larger than $\theta_2$.  Under the null hypothesis of no jump, $\theta_{11}$ is arbitrary, but the parameters $\theta_{12}=0$ and $\theta_2$ disappear; the model is no longer identifiable.  Arbitrary values of $\theta_2$ identify the same distribution for the variable $Y$.
\end{es}
\noindent
Loss of identifiability occurs in areas as diverse as econometrics, reliability theory and survival analysis \citep{PrakasaRao92}, and has been the subject of intensive research.  \cite{Rothenberg71} studied the conditions under which a general stochastic model whose probability law is determined by a finite number of parameters is identifiable.  \cite{PaulinoPereira94} present a systematic and unified description of the aspects of the theory of identifiability.  

As illustrated in Example~\ref{ex:asymptotic}, the classical theory of asymptotic inference heavily relies on quadratic approximation of the log-likelihood function.  Indeed, if the value $\theta_0$ specified by the null hypothesis is unique, we can use \eqref{quadratic} to approximate twice the likelihood ratio function by
\begin{eqnarray}
\nonumber
2\{l(\theta) - l(\theta_0)\} & = & 2\sqrt{n}(\theta-\theta_0)^\top\nu_n(u(\theta_0;Y_i)) \\
\nonumber
            & - & n(\theta-\theta_0)^\top i_1(\theta_0) (\theta-\theta_0) \\
\label{Wquadratic}
            & + & o_p(1),
\end{eqnarray}
for $\theta$ belonging to an $n^{-1/2}$-neighbourhood of $\theta_0$.  Here, $\nu_n(u)=n^{-1/2}\sum_{i=1}^n \{u(\theta; Y_i)-E_{\theta_0}[u(\theta;Y_i)]\}$ is a random process defined for any integrable score function $u(\theta; Y_i)$.  The asymptotic null distribution of $W(\theta_0)$ can be obtained by maximizing this quadratic form in $\theta$.  When the parameter which indexes the true distribution is not unique, various difficulties may arise.  For instance, the maximum likelihood estimator may not converge to any point in the parameter space specified by the null hypothesis.  Or, the Fisher information matrix degenerates.  Typically, the limiting distribution of the likelihood ratio statistics will not be chi-squared.  

In the remainder of the section we will consider two special cases: non-identifiable parameters and singular information matrix.  We will report the main research strains; related contributions are listed in Appendix~B of the Supplementary Material.

\subsection{Non-identifiable parameters}
\label{non-identifiable-parameters}
The general framework for deriving the asymptotic null distribution of the likelihood ratio statistic when some of the parameters are not identifiable under the null hypothesis was developed by \cite{LiuShao03}.  They address the common hypothesis testing problem $H_0: \theta \in \Theta_0$ against $H_1:\theta \in\Theta\setminus\Theta_0$, where $\Theta_0=\{\theta\in \Theta : F_{\theta}=F^0\}$ with $F_\theta$ the distribution function indexed by $\theta$ and $F^0$ the true distribution.  The true distribution is hence unique and $H_0$ is a simple null hypothesis.  However, the set $\Theta_0$ may contain more than one value.  When the true parameter value $\theta^0$ is not unique, the classical quadratic approximation of the likelihood ratio function in an Euclidean neighbourhood of $\theta^0$ no longer holds.  \cite{LiuShao03} by-pass this problem by establishing a general quadratic approximation of the likelihood ratio function, this time in the so-called Hellinger neighbourhood of $F^0$ 
\[ 
\Theta_{\epsilon}=\{\theta\in\Theta \mid 0<H(\theta)\leq \epsilon\},
\] 
where $H^2(\theta)$ is the squared Hellinger distance between $F_\theta$ and $F^0$.  The rationale is that, instead of using the Euclidean distance between two parameter values, $\theta$ and $\theta^0$, closeness between the two models defined by $F_\theta$ and $F^0$ is now measured in terms of a distance, which is valid with or without loss of identifiability of the true distribution $F^0$.  This is closely related to, and indeed generalizes, LeCam’s (1970) local asymptotic theory which we briefly discussed in Section~\ref{local-asymptotics}.  The proof is detailed in Appendix~\ref{identifiability_LiuShao}.  Here, we sketch the main steps in its derivation.

As in \cite{LiuShao03}, we express the likelihood ratio function based on a random sample of $n$ observations $y=(y_1,\ldots,y_n)$
\[
lr(\theta)=\sum_{i=1}^n\log\left\{\lambda_i(\theta)\right\},
\]
in terms of the Radon-Nikodym derivative, $\lambda_i(\theta)=\lambda(\theta;y_i)=dF_\theta/dF^0$, evaluated at $y_i$, for $i=1,\ldots,n$, and recall the definition  
\[
H^2(\theta)=\frac{1}{2}E_{F^0}\left[\left\{\sqrt{\lambda_i(\theta)}-1\right\}^2\right]
\]
of the square Hellinger distance.  As $lr(\theta)$ may diverge to $-\infty$ for some $\theta\in\Theta_\epsilon$, it can be difficult to find a quadratic approximation of the likelihood ratio function with a uniform residual term $o_p(1)$ in $\Theta_\epsilon$.  This is why we rewrite the likelihood ratio statistic as
\begin{equation}
\label{main_W.liushao}
W(H_0) = 2\sup_{\theta\in\Theta\setminus\Theta_0}\{lr(\theta) \vee0 \},
\end{equation}
where $\{a\vee b\} = \max(a,b)$, and maximize $lr(\theta) \vee0$, which generally has a quadratic approximation. 
\cite{LiuShao03} show that such an approximation exists if, for some $\epsilon >0$, the trio $\{S_i(\theta), H(\theta), R_i(\theta)\}$, with $E_{F^0}\left[S_i(\theta)\right]=E_{F^0}\left[R_i(\theta)\right]=0$, satisfies 
\begin{eqnarray*}
h_i(\theta) & = & \sqrt{\lambda_i(\theta)}-1 \\
            & = & H(\theta)S_i(\theta)-H^2(\theta)+H^2(\theta)R_i(\theta), 
\end{eqnarray*}
for all $\theta\in\Theta_\epsilon$, in addition to the generalized differentiable in quadratic mean (GDQM) condition \citep[Definition~2.3]{LiuShao03}.  We may then approximate twice the likelihood ratio function by
\begin{eqnarray}
\label{main_GDQM_exp}
2 \, lr(\theta) & = & 2\sqrt{n}H(\theta)\nu_n(S_i(\theta)) \nonumber \\
       & - & nH^2(\theta)\left\{2+E_{F_n}[S^2_i({\theta})]\right\} + o_p(1),
\end{eqnarray}
where $\nu_n(S_i)=n^{-1/2}\sum_{i=1}^n\left\{S_i(\theta)-E_{F^0}[S_i(\theta)]\right\}$ now is defined in terms of the expectations taken with respect to the empirical distribution function $F_n(\cdot)$ and the true distribution $F^0$.  Expansion~\eqref{main_GDQM_exp} is then used to prove that the distribution of the likelihood ratio statistic \eqref{main_W.liushao} converges to the distribution of the supremum of a squared left-truncated centered Gaussian process with uniformly continuous sample paths.  Though $S_i(\theta)$ and $R_i(\theta)$ may not be unique, they yield the same limiting distribution of the likelihood ratio statistic under suitable conditions. 
In principle, the distribution of the Gaussian process can be approximated by simulation, since its covariance kernel is known.  The most crucial aspect, however, is the derivation of the set which contains the $\mathcal{L}^2$ limits of the generalized score function
\[
\frac{S_i(\theta)}{\sqrt{1+E_{F^0}[S^2_i(\theta)]/2}}
\]
over which the supremum is to be taken.  This needs be worked out on a case by case basis.

The GDQM expansion always exists and reduces to LeCam’s DQM expansion  
\[
h_i(\theta) = (\theta-\theta^0)^\top u(\theta^0;y_i)
\]
if $\theta^0$ is unique.  Furthermore, \cite{LiuShao03} show how \eqref{main_GDQM_exp} is equivalent to \eqref{Wquadratic} by rewriting the latter as
\begin{eqnarray*}
2 \, lr(\theta) & = & 2\sqrt{n}D(\theta^0)\nu_n(S_i^\prime(\theta^0)) \\
                & - & nD^2(\theta^0) + o_p(1),
\end{eqnarray*}
in terms of the squared Pearson-type $\mathcal{L}^2$ distance 
\begin{eqnarray*}
D^2(\theta) & = & E_{\theta^0}\left[ \left\{ \lambda_i(\theta)-1\right\}^2 \right] \\
            & = & (\theta-\theta^0)^\top i_1(\theta^0) (\theta-\theta^0) + o_p(1),
\end{eqnarray*}
if $||\theta-\theta^0|| = O(n^{-1/2})$, and where 
\[
S_i^\prime(\theta) = \left\{ \lambda_i(\theta)-1 \right\} / D(\theta),  
\]
for $\theta\in\Theta \setminus \Theta_0$, defines the generalized score function.

\subsection{Singular information matrix}
\label{singular-information-matrix}

A further case of indeterminate parameter problem is when the expected Fisher information matrix is singular at the true value $\theta^0$ of the parameter.  
\begin{es}[Singular information]
\label{ex:singularity}
Consider a sample of size $n$ from a normal random variable $Y$ with mean $\theta^{q}$, for a given odd integer $q > 0$, and variance $1$.  The information function $i(\theta)=nq^2\theta^{2(q-1)}$ is non singular in an open neighbourhood of $\theta^0$, but vanishes for $\theta^0=0$, which violates Condition~5.  For scalar $\theta$, zero information implies a null score statistic with probability 1.  The left panel of Figure~\ref{fig:4} plots the score functions of three different samples of size $n=10$ for $\theta^0=0$ and $q=3$.  The right panel shows the corresponding normalised log-likelihood functions.  The score function vanishes at the origin and at the maximum likelihoood estimate $\hat\theta=\bar y^{1/q}$.  The log-likelihood function hence admits a global maximum in the neighbourhood of the true parameter value and an inflection at $\theta^0=0$.  Standard techniques to prove consistency of the maximum likelihood estimator and to derive the limiting distribution of the likelihood ratio statistics, such as Expansion~\eqref{Wquadratic}, won’t apply as both $u(\theta_0;y)=i(\theta_0)=0$ at $\theta^0=0$. 
\begin{figure}[tb]
\centering
\includegraphics[width=\linewidth]{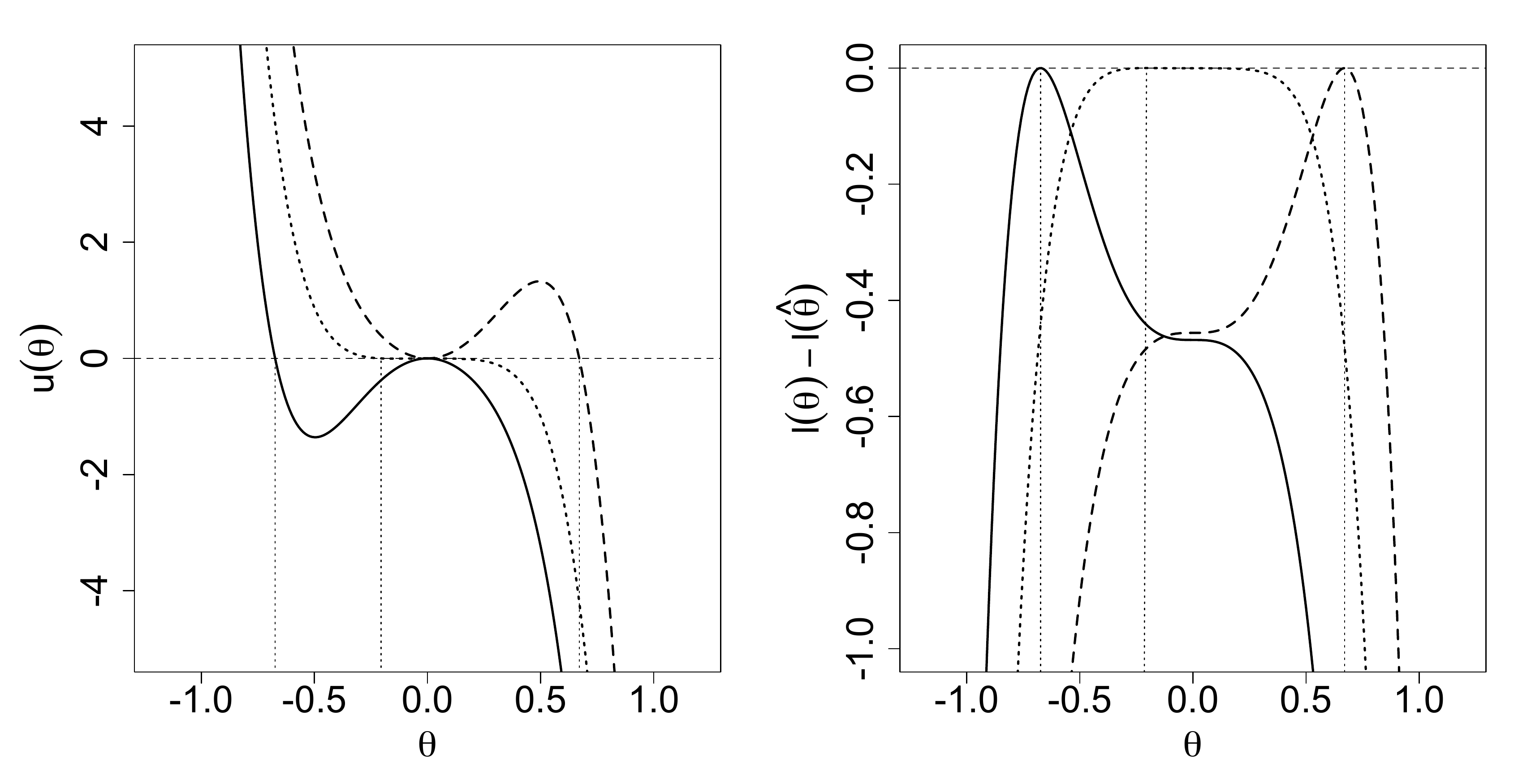}
\caption{\label{fig:4} Example~\ref{ex:singularity}: Singular information.  The three solid curves represent the score functions (left panel) and normalised log-likelihood functions (right panel) for three samples of size $n=10$ for $\theta^0=0$ and $q=3$.}
\end{figure}
\end{es}

Generally speaking, the singularity of the Fisher information matrix prevents the use of the usual second-order expansions of the log-likelihood function.  The, to our knowledge, earliest contribution who addresses this type of problem is \cite{Silvey59}.  The author proposes to modify the curvature of the quadratic approximation of the likelihood ratio by replacing the inverse of the Fisher information matrix with a generalized inverse matrix obtained by imposing suitable constraints on the model parameters.  The however cornerstone contribution to the development of the theory of singular information matrices is \cite{Rotnitzky++00} who derive the asymptotic null distribution of the likelihood ratio statistic for testing the null hypothesis $H_0:\theta=\theta_0$ versus $H_1:\theta\neq\theta_0$, when $\theta$ is a \textit{p}-dimensional parameter of an identifiable parametric model and the information matrix is singular at $\theta_0$ and has rank $p-1$.  Indeed, \cite{Rotnitzky++00} derive a suitable approximation for the likelihood ratio function $l(\theta)-l(\theta_0)$ from a higher-order Taylor expansion.  The theory is developed only for independent and identically distributed random variables, though the authors point out that the same theory may straightforwardly be extended to non-identically distributed observations.  When $\theta$ is scalar, the asymptotic properties of the maximum likelihood estimator and of the likelihood ratio statistic depend on the integer $m_0$, which represents the order of the first derivative of the log-likelihood function which does not vanish at $\theta=\theta_0$; see Theorems~1 and 2 of \cite{Rotnitzky++00}.  If $m_0$ is odd, the distribution of the likelihood ratio converges under the null hypothesis to a $\chi^2_1$ distribution, while for even $m_0$ it converges to a $\bar\chi^2(\omega,1)$ with $\omega=(0.5,0.5)$.  As far as Example~\ref{ex:singularity} goes, $m_0=3$ and $l_3(0;y) = 6n\bar y \neq 0$.  Indeed, the likelihood ratio statistic $W(0) = n \bar Y^2$ distributes exactly as a chi-squared distribution with one degree of freedom.  Extensions of these results when the parameter $\theta$ is $p$-dimensional are also provided.  These are generally based on suitable re-parametrizations of the model which remove the specific causes of the singularity, but are difficult to generalize as they are ad-hoc solutions.  Further contributions who propose to penalise the likelihood function so as to guarantee the consistency and normality of the maximum likelihood estimator are mentioned in Appendix~B of the Supplementary Material.


\section{Finite mixture models}
\label{finite-mixture-models}

\subsection{Background}
Finite mixtures deserve special attention, because of their widespread use in statistical practice, but also because of the methodological challenges posed by the derivation of their asymptotic properties.  They probably represent the best-studied indeterminate parameter problem, though we may also treat them as a boundary case.  This is best illustrated by the two-component mixture model 
\begin{equation}
\label{2mixture}
f(y;\eta)=(1-\pi) f_1(y;\theta_1) + \pi f_2(y;\theta_2), 
\end{equation}
with $\eta=(\pi,\theta_1,\theta_2)$. 
Here, the probability density or mass functions $f_1(y;\theta_1)$ and $f_2(y;\theta_2)$, with $\theta_1\in\Theta_1\subseteq{\mathbb R}^{p_1}$ and $\theta_2\in\Theta_2\subseteq{\mathbb R}^{p_2}$, represent the mixture components, while $0\leq\pi\leq 1$ is the mixing probability.  The null hypothesis of homogeneity can be written in different ways.  We may set $\pi=0$, which corresponds to $H_0: f^0 = f_1(y;\theta_1)$, where $f^0$ represents the true unknown distribution, or alternatively, $\pi=1$ and $H_0: f^0 = f_2(y;\theta_2)$.  If the two components, $f_1(y;\theta_1)$ and $f_2(y;\theta_2)$, are known, then the limiting distribution is a $\bar\chi^2(\omega,1)$ with $\omega=(0.5,0.5)$ \citep[p.~75]{Lindsay95}.  Otherwise, for $f_1(y;\theta) = f_2(y;\theta)$ a third possibility arises: in this case homogeneity assumes that $H_0: \theta_1 = \theta_2$.  Whatever choice is made, some model parameters, that is, $\theta_2$ and $\theta_1$, respectively, in the first two cases and $\pi$ in the third, vanish under the null hypothesis.  This contradicts classical likelihood theory, where the parameter which characterises the true distribution is typically assumed to be a unique point $\theta^0$ in the open subset $\Theta\subseteq{\mathbb R}^p$.  Under this scenario, the asymptotic distribution of the likelihood ratio statistic does not follow the commonly believed chi-squared distribution.  Indeed, in many cases the finite sample distribution converges to the supremum of a Gaussian process. 

The remainder of the section outlines the mainstream contributions for this class of models, with special emphasis on homogeneity testing using the likelihood ratio.  This represents a subset of all available and related contributions on finite mixtures, whose treatment would easily fit a book-length account.  Two general references for finite mixture models are the monographs by \cite{Lindsay95} and \cite{McLachlanPeel00}.  The large-sample properties of a number of classical and recent likelihood-based test statistics for assessing the number of components of a finite mixture model are reviewed in \cite{Chen17}.  Further related work is listed in Appendix~B of the Supplementary Material.

\subsection{Testing for homogeneity}
\label{testing-for-homogeneity}

\subsubsection{General theory}
The first discussion of asymptotic theory for testing homogeneity of model~(\ref{2mixture}) when all parameters are unknown was provided by \cite{GhoshSen85} who characterise the limiting distribution of the likelihood ratio statistic under the assumption that $\Theta_2$ is a closed bounded interval of ${\mathbb R}$, while $\Theta_1\subseteq{\mathbb R}^{p_1}$, $p_1\geq 1$.  There is an additional major difficulty in dealing with finite mixture models: though the mixture itself may be identifiable, the parameters $\pi$, $\theta_1$ and $\theta_2$ may not be.  Indeed, if $f_1(y;\theta) = f_2(y;\theta)$ in \eqref{2mixture}, the equality 
\begin{equation*}
f(y;\eta)=f(y;\eta'),
\end{equation*}
holds for $\eta=\eta^\prime=(\pi^\prime,\theta_1^\prime,\theta_2^\prime)$, but also for
$(1-\pi^\prime,\theta_2^\prime,\theta_1^\prime)$.  That is, under the alternative hypothesis there is a second set of parameters which gives rise to the same distribution, while under the null hypothesis of homogeneity the model is represented by the three curves $\pi = 1$, $\pi =0$ and $\theta_1=\theta_2$.  Choosing an identifiable parametrisation doesn't bring any improvement as the density is then no longer differentiable.  \cite{GhoshSen85} by-pass this difficulty in two ways: by requiring either strong identifiability of the finite mixture or then by imposing a separation condition on the parameters.

Strong identifiability holds if the equality $f(y;\eta)=f(y;\eta^\prime)$ implies that $\pi=\pi^\prime$, $\theta_1=\theta_1^\prime$ and $\theta_2=\theta_2^\prime$.  The distribution of the likelihood ratio statistic for testing $H_0: \pi = 0$ then converges to the distribution of $T^2I_{\{T>0\}}$, where $T=\sup_{\theta_2}\{Z(\theta_2)\}$ and $Z(\theta_2)$ is a zero-mean Gaussian process on $\Theta_2$ whose covariance function depends on the true value of the parameters under the null hypothesis \citep[Theorem~2.1]{GhoshSen85}.  This results from proceeding in two steps.  We first approximate the log-likelihood function by a quadratic expansion with respect to $\pi$ and $\theta_1$ which, under the null hypothesis, converges to the square of a Gaussian random process indexed by the non-identifiable parameter $\theta_2$.  The supremum of this process with respect to $\theta_2$ is then taken.  The sketch of this proof is given in Appendix~\ref{finite_mixture_GhoshSen85}.

A similar result holds if the finite mixture is not strongly identifiable, such as when $f_1(y;\theta) = f_2(y;\theta)$ in (\ref{2mixture}).  In this case, a separation condition between $\theta_1$ and $\theta_2$ of the form $||\theta_1-\theta_2|| \geq \epsilon$ for a fixed quantity $\epsilon > 0$ needs be imposed, so that $H_0$ is described by either $\pi=0$ or $\pi=1$ \citep[\S 5]{GhoshSen85}.  The proof outlined in Appendix~\ref{finite_mixture_GhoshSen85} still applies with the exception that now the non-identifiable parameter $\theta_2$ varies in a subset of $\Theta_2$ which depends on the given $\epsilon$.

\subsubsection{Gaussian mixtures}
\label{gaussian-mixture}
Theoretical results are particularly generous if the two-component model is a normal mixture, which is justified by the widespread use of the normal distribution in a wide variety of situations.  Milestone contributions are \cite{Goffinet++92} and \cite{ChenChen01, ChenChen03}, who consider a two-component mixture of normal densities $\phi(\cdot;\mu,\sigma)$ of mean $\mu$ and variance $\sigma^2$.  This class of models deserves special attention because of the technical challenges posed by some undesirable properties of the Gaussian distribution.  As discussed in \cite{ChenLi09}, the expected Fisher information for the mixing proportion is not finite unless a bound is placed on the variance of the corresponding component density.  Furthermore, the derivatives of the log-density become linearly dependent since $\partial^2\phi(x;\mu,\sigma)/\partial\mu^2=2\partial\phi(x;\mu,\sigma)/\partial\sigma^2$.  Last, but not least, the log-likelihood function is unbounded in case of heterogeneous components, and the maximum likelihood estimator may not exist \citep{Hartigan85}.  These deficiencies generally invalidate the standard quadratic approximation of the likelihood ratio function, which needs be expanded further.  This is why most contributions on normal mixtures assume homogeneous variances.  

Several solutions have been proposed to overcome these shortcomings; in addition to \cite{ChenLi09}, see also \cite{KasaharaShimotsu15} and \cite{Wichitchan++19}.  The finite sample distribution of the likelihood ratio statistic often converges to the distribution of
\begin{equation}
\label{chenchen}
\big\{\!\! \sup_{|t|\leq M}Z(t) \big\}^2+ W,
\end{equation}
where $Z(t)$, $t\in[-M,M]$, is a Gaussian process and $W$ is an independent chi-squared random variable with suitable degrees of freedom.  The Gaussian process $Z(t)$ has zero mean and known covariance function.  As mentioned in Section~\ref{regularity-conditions}, the compactness of the parameter space is a necessary condition to avoid that the distribution of the likelihood ratio statistic diverges to infinity.  This was already proved by \cite{Hartigan85} and is an immediate implication of the fact that $\{\sup_{|t|\leq M}Z(t)\}^2$ tends in probability to infinity if $M\rightarrow\infty$.  The proofs of the Theorems in \cite{ChenChen01,ChenChen03} essentially are suitable adaptations of the prototype derivation for finite mixture models reported in Appendix~\ref{finite_mixture_GhoshSen85}.  All passages are detailed in the original contributions to which we refer the interested reader.  

As in most cases the asymptotic distribution of the likelihood ratio statistic is related to a Gaussian random field, the computation of percentile points becomes tricky or impossible.  That is why other tests or methods have been proposed.  Reviewing all these would go beyond the scope of the paper.  Let us mention, here, the most fruitful research strained initiated by \cite{Li++09} who propose an EM-test for homogeneity, which \cite{ChenLi09} decline in the case of a two-component Gaussian mixture.  A most recent treatment is \cite{Chauveau++18}.

\subsection{Alternative approaches}
\label{relaxing-the-conditions}
Several authors have addressed how to remove the separation condition of \cite{GhoshSen85}.  Three lines of research emerged: reparametrization of the probability density or mass function, penalization of the likelihood to ensure identifiability and simulation.  

Reparametrization does not change the model which is why the limiting distribution, generally, remains the supremum of a Gaussian process.  The first contribution, to our knowledge, who uses reparametrization is \cite{ChernoffLander95} who heuristically study several versions of the two-component binomial mixture model.  Formal proofs and extensions to finite mixtures with contaminated densities are provided in \cite{LemdaniPons97,LemdaniPons99}, while \cite{Ciuperca02} considers the case of translated mixture components.  This latter contribution has the further merit of highlighting how Condition~3 of Section~\ref{regularity-conditions} is necessary, but not sufficient.  Indeed, the limiting distribution of the likelihood ratio statistic converges to a fifty-fifty mixture of a point mass at zero and of a distribution which diverges in probability to $+\infty$, and this despite the fact that all parameters are assumed to belong to a compact set.  The unboundedness behaviour of the likelihood ratio of \cite{Ciuperca02} can be explained by means of the theory of ``locally conic\rq\rq\ reparametrizations proposed by \cite{Dacunha-CastelleGassiat97,Dacunha-CastelleGassiat99}. 

A rather different route is taken in \cite{Chen++01} who suggest to penalise the log-likelihood function of model~\eqref{2mixture} with $f_1(y;\theta) = f_2(y;\theta)$,
\begin{equation}
\label{loglik_m}
l(\pi,\theta;y)+ c\log\{4\pi(1-\pi)\},
\end{equation}
where the degree of penalisation is controlled by the constant term $c$.  As the authors point out, the penalisation term can be justified from the Bayesian perspective, as a prior on the mixing proportion $\pi$.  It furthermore guarantees that the maximum likelihood estimate of the mixing proportion $0<\hat\pi<1$ will not fall on the boundary of the parameter space and that the maximum likelihood estimators of all parameters are consistent under the null hypothesis $\pi=0$.  Provided Conditions~1--5 of their paper hold, the distribution of the modified likelihood ratio statistic derived from (\ref{loglik_m}) converges to a $\bar\chi^2(\omega,1)$ distribution with $\omega=(0.5,0.5)$; see also Example~\ref{ex:vonMises}.  Indeed, regularising the likelihood function does change the problem at hand and the limiting distributions no longer is the supremum of a squared truncated Gaussian random process. 

The third route to investigate the asymptotic null distribution of the likelihood ratio statistic for finite mixture models is by simulation.  \cite{Thode++88} consider testing the hypothesis that the sample comes from a normal random variable with unknown mean and unknown variance against the alternative that the sample comes from a two-component Gaussian mixture with unequal means and common variances.  All model parameters are assumed to be unknown.  Their extensive numerical investigation shows that the distribution of the likelihood ratio statistic converges very slowly to a limiting distribution, if any exists, and is rather unstable even for sample sizes as large as $n=1,000$.  For very large sample sizes, the empirical distributions rather closely agree with the commonly assumed $\chi^2_2$, though this may be too liberal for small to moderate $n$.  \cite{Bohning++94} investigate numerically the asymptotic properties of the likelihood ratio statistic for testing homogeneity in the two-component mixture model~(\ref{2mixture}) when the component distributions $f_k(y;\theta_k)$, $k=1,2$ are binomial, Poisson, exponential or Gaussian with known common variance.  \cite{Lo08} shows that the commonly used $\chi^2$ approximation for testing the null hypothesis of a homoscedastic normal mixture against the alternative that the data arise from a heteroscedastic model is reasonable only for samples as large as $n=2,000$ and component distributions that are well separated under the alternative.  Furthermore, the restrictions of \cite{Hathaway85} need be imposed to ensure that the likelihood is bounded and to rule out spurious maxima under the alternative.  Otherwise, the author suggests use of parametric resampling.  Very recently, \cite{CongYao21} study the behaviour of the likelihood ratio statistics for multivariate normal mixtures.  They recommend using parametric boostrap resampling as, similarly to \cite{Lo08}.


\section{Change-point problems}
\label{change-point-problems}

\subsection{Definition}
\label{change-point-problems-background}
A change-point problem arises whenever the regime of random events suddenly changes.  A modification in the data generating process generally implies that the log-likelihood function is no longer differentiable with respect to some values of the parameter.  This typically leads to the failure of Condition~4 of Section~\ref{regularity-conditions}.  The most basic change-point problem tries and identifies patterns in a random sequence.  For instance, given $n$ independent observations $y_1,\ldots, y_n$, listed in the order they occurred, \cite{Page55, Page57} considered the problem of verifying whether these were generated by a random variable with distribution function $F(y;\theta)$ against the alternative that only the first $\tau$, $0\leq\tau<n$, observations are generated from $F(y;\theta)$ while the remaining $n-\tau$ come from $F(y;\theta^\prime)$ with $\theta\neq\theta^\prime$ and $\tau$ unknown.  

Since Page’s pioneering papers, change-point problems have been the subject of intensive research owing to their pervasiveness in all major domains of application.  Summarizing all contributions can easily fill in book-length accounts.  A first annotated bibliography of change-point problems is \cite{Shaban80}.  \cite{KrishnaiahMiao88} give an overview of change-point estimation up to their time of writing; \cite{CsorgoHorvath97} focus their review monograph on limit theorems for change-point analysis.  \cite{KhodadadiAsgharian08} is a more than 200 pages length annotated bibliography of change-point problems in regression.  \cite{Lee10} summarizes the most recent literature and gives a comprehensive bibliography for five major types of change-point problems.  A book-length account of change-point problems with examples from medicine, genetics and finance is \cite{ChenGupta12}.  The discussion papers of \cite{HorvathRice14a,HorvathRice14b} mention, in addition to classical methods, also modern lines of research in functional data and high dimensions.  Research on change-point analysis has seen a revival over the last decade and a half, especially as far as the detection of multiple changes goes \citep{Niu++16,YauZhao16,DetteGosmann20}.  The proposed inferential solutions range from parametric to nonparametric techniques and include frequentist and Bayesian approaches.  Most recently, \cite{Sofronov++20} edited a special issue of the journal \textit{Statistical Papers} dedicated to change-point detection.

Generally speaking, two questions are of interest in change-point analysis: identifying the potential number of changes, and, once identified, estimating where these occur, together with further quantities of interest such as the size of the change.  In the remainder of this section we will focus our attention on the first problem, that is, the identification of a change by means of the likelihood ratio statistic.  As highlighted by \cite{ChenGupta12}, the majority of reference models which have been proposed for change-point detection assume normality of the observations.  These will be treated extensively in Sections~\ref{shifts-in-location-and-dispersion} and \ref{change-point-detection-in-regression} with special emphasis on regression type problems.  In particular, Section~\ref{shifts-in-location-and-dispersion} addresses the issue of detecting possible shifts in the location and/or the scale of the distribution.  Sections~\ref{change-point-detection-in-regression} extends the treatment to linear regression and piecewise linear models.  Given the breadth of the available solutions, each section contains a selection of contributions which illustrate the main currents of research.  A further deeply explored class of models are continuous exponential families \citep{Worsley86}, though some results are available for the binomial \citep{Worsley83} and Poisson cases.  Further contributions are listed in Appendix~B of the Supplementary Material.

\subsection{Shifts in location and scale}
\label{shifts-in-location-and-dispersion}
The reference model for testing a change in the mean value of a random variable can generally be written as
\begin{equation}
\label{regression}
y_i = \eta_i + \varepsilon_i, \quad i=1,\ldots,n , 
\end{equation}
where the $\varepsilon_i$'s are independent zero-mean random errors.  All observations are considered in the order they appear, an assumption which will hold for the whole section.  The function $\eta_i$ may change $K$ times,
\begin{eqnarray}
\label{cp-model}
\eta_i & = & \mu_1, \quad 0 < i \leq \tau_1, \\
\nonumber
       & = & \mu_2, \quad \tau_1 < i \leq \tau_2, \\
\nonumber
       & \vdots & \\
\nonumber
       & = & \mu_{K+1}, \quad \tau_K < i \leq n, 
\end{eqnarray}
where the change-points $\tau_k$ can only assume integer values.  Unless differently stated, both the $K+1$ different mean values $\mu_k$ and the $K$ change-points $\tau_k$ are supposed to be unknown, though the very early contributions focus on the simpler setting where one or both pieces of information are given. 

Assuming $K=1$, \cite{Hawkins77} considers testing the null hypothesis of no change in the mean $\eta_i$ when the $\varepsilon_i \sim N(0, \sigma^2)$ are centered normal variables with constant variance $\sigma^2>0$, that is,
$$H_0: Y_i\sim N(\mu,\sigma^2), \quad i=1,\ldots,n,$$
against the alternative that there exists a $0<\tau<n$ at which the unknown mean switches from $\mu$ to $\mu^\prime\neq\mu$.  The variance $\sigma^2$ is assumed to be known and we set it to one without loss of generality.  This is a non-standard problem because the change-point appears only under the alternative hypothesis, but not under the null.  The corresponding likelihood ratio statistic is a function of 
\[ 
W_\tau = \tau (\bar Y_\tau - \bar Y)^2 + (n - \tau)(\bar Y_{n-\tau}-\bar Y)^2,
\]
where $\bar Y_\tau$ and $\bar Y_{n-\tau}$ are the partial means, computed using the first $\tau$ and the last $n-\tau$ observations.  Since the change-point $\tau$ is unknown, the likelihood ratio 
\[
W = W_{\tau^*}=\max_{1\leq \tau< n} W_\tau  
\]
maximises $W_\tau$ over all possible values of $\tau$, and is usually refereed to a \lq\lq maximally selected likelihood ratio''.  To derive its exact null distribution, \cite{Hawkins77} re-expresses $W_\tau$ as  
\[
W_\tau=T_\tau^2,
\]
where
\begin{equation}
\label{eq:t}
T_\tau=\sqrt{\frac{n}{\tau(n-\tau)}}\sum_{i=1}^\tau(Y_i-\bar{Y})
\end{equation}
has standard normal distribution.  It follows that the finite-sample distribution of 
\begin{equation}
\label{U_statistics}
U = \sqrt{W_{\tau^*}}=\max_{1\leq \tau< n}|T_\tau|
\end{equation}
agrees with the distribution of the maximum absolute value attained by a Gaussian process in discrete time having zero mean, unit variance and autocorrelation function given by Expression~(3.2) of \cite{Hawkins77}.  In particular, the null distribution of $U$ has density function 
\begin{equation}
\label{u_density}
f_{U}(u)=2\phi(u)\sum_{\tau=1}^{n-1}g_{\tau}(u)g_{n-\tau}(u),
\end{equation}
where $\phi(u)$ is the density of the standard normal, $g_1(u)=1$ for $u\geq 0$ and $g_\tau(u)$ is a recursive function such that  
\begin{equation}
\label{eq:g}
g_{\tau}(u)={\rm Pr}(|T_i|<u, i=1,\ldots, \tau-1\mid |T_\tau|=u).
\end{equation}
The sketch of the proof of \eqref{u_density} is given in Appendix~\ref{gaussian_mean_shift_Hawkins77}.  \citet[Theorem~2.1]{YaoDavis86} show that a suitably normalized version of $U$ converges, however slowly \citep{Jaruskova97,CsorgoHorvath97}, under $H_0$ to the double exponential, or Gumbel, distribution, which provides approximate quantiles.  See also \cite{Irvine86}.

The finite-sample null distribution of the likelihood ratio statistic for unknown $\sigma^2$ was worked out by \cite{Worsley79}.  Further generalizations, such as to the multivariate case and/or to account for a possible change in the scale of the distribution, can be found in the book length account of \citet[\S\S 2.2--2.3 and 3.2--3.3]{ChenGupta12}.  See also the selection of references given in Appendix~B of the Supplementary Material.

\subsection{Change-point detection in regression}
\label{change-point-detection-in-regression}
A further extension of Model~(\ref{cp-model}) with respect to location,
\begin{eqnarray}
\label{time-varying}
\eta_i & = & \alpha_1 + \beta_1 x_i, \quad 0 < i \leq \tau_1, \\
\nonumber
       & = & \alpha_2 + \beta_2 x_i, \quad \tau_1 < i \leq \tau_2, \\
\nonumber
       & \vdots & \\
\nonumber
       & = & \alpha_{K+1} + \beta_{K+1} x_i, \quad \tau_K < i \leq n,
\end{eqnarray}
is used for change-point detection in simple linear regression.  The early contributions by \cite{Quandt58,Quandt60} derive the likelihood ratio statistic under the null hypothesis of no switch against the alternative that the model possibly obeys two separate regimes under the assumption of independent and zero-mean normal error terms $\varepsilon_i$.  Under the alternative hypothesis, the variance is furthermore allowed to switch from $\sigma_1^2$ to $\sigma_2^2$ at instant $\tau$, when the linear predictor $\eta_i$ undergoes a structural change.  The likelihood ratio statistic is
\begin{equation}
\nonumber
W = \max_{3\leq \tau\leq n-3} W_\tau,
\end{equation}
with
$$
W_\tau=-2\log\left(\frac{\hat\sigma_1^{2\tau}\; \hat\sigma_{2}^{2(n-\tau)}}{\hat\sigma^{2n}}\right),
$$
is a function of the least squares estimators $\hat{\sigma}^2_{1}$ and $\hat{\sigma}^2_{2}$ of $\sigma^2_1$ and $\sigma^2_2$, respectively, computed using the corresponding subsets of observations, and of the MLE $\hat{\sigma}^2$ of the common variance $\sigma^2 = \sigma_1^2 = \sigma_2^2$ based upon the entire sample. \cite{Quandt58} initially conjectured that the asymptotic distribution of $W$ may be $\chi^2_4$ under the null hypothesis of no change.  However, the numerical investigation he reported in a later publication for the three sample sizes $n=20, 40, 60$ \citep[Table~3]{Quandt60} revealed that the finite-sample distribution depends on the number of observations $n$. 

Change-point detection in simple linear regression using the likelihood ratio is also the subject of \cite{KimSiegmund89}.  These authors consider two situations: where only the intercept is allowed to change and where both, the intercept and the slope change while the variance remains constant.  Again, the Brownian Bridge process is central to the derivation of the corresponding limiting distributions as in \cite{YaoDavis86}.  Approximations for the corresponding tail probabilities are given by \cite{KimSiegmund89} under reasonably general assumptions. 

Model~\eqref{time-varying} can be extended to account for changes in the covariates, which is known as piecewise linear regression.  This type of models are very popular in a large number of disciplines, including among others environmental sciences \nocite{PiegorschBailer97,Muggeo08a} (Piegorsch and Bailer, 1997, Section~2.2; Muggeo, 2008a), medical sciences \citep{SmithCook80,Muggeo++14}, epidemiology \citep{Ulm91} and econometrics \citep{Zeileis++06}.  A review of likelihood ratio testing for piecewise linear regression up to his time of writing can be found in \cite{Bhattacharya94}.  See also the annotate bibliography for a selection of related contributions.


\section{Beyond parametric inference}
\label{beyond-parametric-inference}
This section reviews cases of interest which do not fit into the previously mentioned three broad model classes, but still fall under the big umbrella of non-standard problems.  In particular we will focus on shape constrained inference, a genre of nonparametric problem which leads to highly nonregular models.  

As brought to our attention by an anonymous Referee, the asymptotic theory of semiparametric and nonparametric inference has interesting analogues to the classical parametric likelihood theory reviewed in Section~\ref{likelihood-asymptotics}.  Indeed, the parameter space of a semiparametric model is an infinite-dimensional metric space.  This makes the model non-standard as we typically consider a real parameter of interest in the presence of an infinitely large nuisance parameter.  Despite this departure from regularity, the likelihood ratio statistic still behaves as we would expect it.  \cite{MurphyVaart97,MurphyVaart00}, for instance, show that the corresponding limiting distribution is chi-squared also when we profile out the infinite-dimensional nuisance parameter.  The behavior of this likelihood ratio statistics under local alternative hypotheses is studied in \cite{Banerjee05}.  The classical approximations of Section~\ref{likelihood-asymptotics} also hold for the asymptotic theory of empirical likelihood \citep{Owen90,Owen91}; see \cite{ChenKeilegom09} for a review.  These results are quite remarkable given that the underlying distributional assumptions are much less strict.

An area of research which has received much attention in the last decade is nonparametric inference under shape constraints \citep{StatScience2018}.  Shape constraints originate as a natural modelling assumption and lead to highly nonregular models.  As highlighted by \cite{GroeneboomJongbloed18}, the probability density or mass functions of many of the widely used parametric models satisfy shape constraints.  For example, the exponential density is decreasing, the Gaussian density is unimodal, while the Gamma density can be both, depending on whether its shape parameter is smaller or larger than one.  Estimation under shape constraints leads to an M-estimation problem where the parameter vector typically has the same length as the sample size and is constrained to lie in a convex cone.  Non-regularity arises since the M-estimator typically falls on the face of the cone.  As for boundary problems, convex geometry is an essential tool to treat shape constrained problems.  

The field of shape constraint problems originated from `monotone’ estimation problems, where functions are estimated under the condition that they are monotone.  The maximum likelihood estimator converges typically at the rate $n^{-1/3}$ if reasonable conditions hold, that is, at a slower pace than the $n^{-1/2}$ rate attained by regular problems.  Moreover, the maximum likelihood estimator has a non-standard limiting distribution known as Chernoff’s distribution \citep{GroeneboomWellener01}.  A considerable body of work has studied the asymptotic properties of the nonparametric likelihood ratio statistic under monotonicity.  In particular, \cite{BanerjeeWellner01} initiated the research strain of testing whether a monotone function $\psi$ assumes the particular value $\psi(t_0)=\psi_0$ at a fixed point $t_0$.  An extension to regression is given by \cite{Banerjee07}, who assumes that the conditional distribution $p(y,\theta(x))$, of the response variable $Y$ given the covariate $X=x$, belongs to a regular parametric model, where the parameter $\theta$, or part of it, is specified by a monotone function $\theta(x)\in\Theta$ of $x$.  Other types of shape constraint problems have emerged in the meantime which entail concavity or convexity and uni-modality of the functions to be estimated; see Appendix~B of the Supplementary Material.  

Shape constraints arise also in many high-dimensional problems, which opens frontiers for research in nonregular settings; see for example \citet{Bellec18}.  Most recently, \cite{DossWellner19} showed that the likelihood ratio statistic is asymptotically pivotal if the density is log-concave.  The class of log-concave densities has many attractive properties from a statistical viewpoint; an account of the key aspects is given in \cite{Samworth18}. Non-standard limiting  distributions characterize shape constrained inferential  problems.  Generally, the likelihood ratio statistic converges to a limiting distribution which can be described by a functional of a standard Brownian motion plus a quadratic drift.  In addition, the limiting distribution is asymptotically pivotal, that is, it doesn’t depend on the nuisance parameters, as happens for the common $\chi^2$ distribution of regular parametric problems.  Recent work on high-dimensional asymptotics of likelihood ratio tests under convexity constraints is discussed in \cite{Han++22}.


\section{Computational aspects and software}
\label{software}

Deriving the asymptotic distribution of the likelihood ratio statistic under non-standard conditions is generally a cumbersome task.  In some cases the limiting distribution is well defined and usable, as for instance when it boils down to a chi or chi-bar squared distribution.  Quite often, however, the analytical approximation is intractable, so as when we have to determine the percentiles of a Gaussian random field.  This fact has motivated the development of alternative test statistics whose null distribution presents itself in a more manageable form; see, for instance, the contributions mentioned in Section~\ref{relaxing-the-conditions}.  Or, we may rely upon simulation, using Monte Carlo or the bootstrap, as mentioned in passing in Sections~\ref{boundary_null_variance}, \ref{non-identifiable-parameters}, \ref{relaxing-the-conditions} and \ref{shifts-in-location-and-dispersion}.  Bootstrapping, in particular, allows us to recover the finite-sample null distribution of the test statistic very naturally provided that the bootstrap resamples from a consistently estimated density \citep{Titterington90, FengMcCulloch96}.  In general terms, this requires that the maximum likelihood estimator converges to the possibly non-identifiable subset of the parameter space to which the true parameter belongs to.  An example of inconsistency of the bootstrap when a parameter is on the boundary of the parameter space together with further counterexamples that are already in the literature are provided in \cite{Andrews00}.  Bootstrap likelihood ratio tests for finite-mixture models are reviewed in \cite{FengMcCulloch96}.  A most recent application for boundary points is \cite{Cavaliere++22}, while \cite{Kirch08}, \cite{HusskovaKirch12}, and the very recent papers by \cite{ChenCabrera20} and \cite{YuChen22} boostrap the critical values of change-point tests.  Permutation is also used to derive the critical values for tests statistics in change-point analysis; see for example \cite{KirchSteinebach06} and reference therein.

A compromise between analytical approximation and simulation is the hybrid approach described in \citet[Section~7.7]{Brazzale++07} where parts of the analytical approximation are obtained by simulation.  However, simulation becomes useless if the limiting distribution diverges to infinity, as already mentioned in Example~\ref{ex:vonMises}.  A non exhaustive list of examples is provided in paragraphs~\ref{gaussian-mixture} and \ref{relaxing-the-conditions} of Appendix~B of the Supplementary Material.  Substantive applications in which the approximations have been found useful and details of how to implement the methods in standard computing packages are generally missing.  Reviewing all software contributions which implement likelihood ratio based inference for nonregular problems in a more or less formalized way is beyond the scope of this paper.  In the following we try and give a selected list of packages for the numerical computing environment \textsf{R} \citep{RCoreTeam20}.  We will again group them into the three broad classes reviewed in the previous Sections~\ref{boundary-problems}--\ref{change-point-problems}, that is, boundary problems, mixture models and change-point problems. 
 
Crainiceanu and Ruppert’s (2004) \nocite{CrainiceanuRuppert04} proposal, which tests  for a null variance component, is implemented in the \texttt{RLRsim} package by \cite{Scheipl++08}.  We furthermore mention the \texttt{varTestnlme} package by \cite{BaeyKuhn19} and the \texttt{lmeVarComp} package by \cite{Zhang18}.  The first again tests for null variance components in linear and non linear mixed effects model, while the second implements the method proposed by \cite{Zhang++16} for testing additivity in nonparametric regression models.  

An account of some early software implementations to handle mixture models can be found in \cite{Haughton97}, in the Appendix of \cite{McLachlanPeel00} and also in the Software section of the recent review paper by \cite{McLachlan++19}.  A most recent implementation for use in astrostatistics is the \texttt{TOHM} package by \cite{AlgerivanDyk20} which implements a computationally efficient approximation of the likelihood ratio statistic for a multimensional two-component finite-mixture model.  The package is also available for the Python programming language.  The code provided by \cite{Chauveau++18} for testing a two-component Gaussian mixture versus the null hypothesis of homogeneity using the EM test is available through the \texttt{MixtureInf} package by \cite{Li++16}.  Maximum likelihood estimation in finite mixture models based on the EM algorithm is furthermore addressed in the \texttt{mixR} package by \cite{Yu18}, which also considers different information criteria and bootstrap resampling.  The \texttt{clustBootstrapLRT} function of the \texttt{mclust} package by \cite{Scrucca++16} also implements bootstrap inference for the likelihood ratio to test the number of mixture components.  A further implementation of the likelihood ratio test for mixture models is the \texttt{mixtools} package by \cite{Benaglia++09}.  All \textsf{R} packages linked to finite mixture models are listed on the CRAN Task View webpage for Cluster Analysis \& Finite Mixture Models\footnote{\url{http://cran.r-project.org/web/views/Cluster.html}}.

The \texttt{changepoint} package by \cite{KillickEckley14} considers a variety of test statistics for detecting change-points among which the likelihood ratio.  The \texttt{strucchange} package by \cite{Zeileis++02} provides methods for detecting changes in linear regression models.  We may furthermore mention the \texttt{segmented} package by \cite{Muggeo08b} for change-point detection in piecewise linear models, the \texttt{bcp} package by \cite{ErdmanEmerson07} for Bayesian analysis of a single change in univariate time series and the \texttt{CPsurv} package by \cite{Brazzale++19} for nonparametric change-point estimation in survival data.


\section{Discussion}
\label{discussion}
Non-regularity can arise in many different ways, though all entail the failure of one, at times even two, regularity conditions.  Many problems can be dealt with straightforwardly; other require sophisticated tools including limit theorems and extreme value theory for random fields.  The wealth of contributions, which has been produced during the last 70 years, synthesized in Figure~\ref{timeline}, testifies that the interest in this type of problems has not faded since they made their entrance back in the early 50’s.  The best-studied nonregular case are boundary problems.  Common examples of application are testing for a zero variance component in mixed effect models and constrained one-sided tests.  The limiting distribution of the likelihood ratio is generally a chi-bar squared distribution with a number of components and mixing weights that depend on the number of parameters which fall on the boundary.  This is also the only type of problem for which higher order results are available.

Indeterminate parameter problems are far more heterogenous.  Apart from finite mixtures, the remaining cases can be put under the two umbrellas of non-identifiable parameters and singular information matrix.  The methodological difficulties increase as the limiting distributions depend on the parametric family and on the unknown parameters.  If $\theta$ is scalar and we want to test homogeneity against a two-component mixture, the distribution of the likelihood ratio converges to the distribution of the supremum of a Gaussian process.  For a larger number of mixture components and/or multidimensional $\theta$, this becomes the distribution of the supremum of a Gaussian random field.  In these cases, simulation-based approaches are often needed to obtain the required tail probabilities.  Moreover, constraints must be imposed to guarantee identifiability of the mixture parameters.  As outlined by \cite{Garel07}, these may act on the parameter space, by bounding it or imposing suitable separation conditions among the parameters, or on the alternative hypotheses which must be contiguous.  A further possibility is to penalize the likelihood function so that the limiting distribution of the corresponding modified likelihood ratio statistic is chi-squared or well approximated by a chi-bar squared distribution.

Change-point problems range from the simple situation of detecting an alteration in the regime of a random sequence to identifying a structural break in multiple linear regression with possibly correlated errors.  Although in the latter case the change-point can assume any value, in the first situation it must lie in a discrete set.  In addition, there is a tie between change-point problems and indeterminate parameter problems whenever setting one of the components of the model to a particular value, can make other components, or parts of it, disappear, as shown in Example~\ref{jump-regression}.  Generally, the limiting distribution of the likelihood ratio statistic for detecting a change either converges to a Gaussian processes or can be adjusted to converge to a Gumbel type distribution \citep{HorvathRice14a}.  This technique was first applied by \cite{DarlingErdos56} to derive the limiting distribution of the maximum of independent random variables, and has further been extended to depedent data; see \cite{AueHorvath13} for a review.  Approximate critical values of the test statistics can be obtained from Bonferroni’s inequality, by using asymptotic arguments or simulation.  In some situations the likelihood ratio statistic for the unknown change-point is unbounded. 

From the more practical point of view, use of the asymptotic distribution of the likelihood ratio statistic loses its appeal once it goes beyond the common $\chi^2$ distribution.  As a result, simulation-based tests that circumvent the asymptotic theory are often used.  Indeed, simulation may nowadays be used to establish the desired empirical distributions of the estimators and to compute approximations for $p$-values obtained from Wald-type statistics.  For the most intricate situations, the authors suggest to use resampling-based techniques, such as parametric and nonparametric bootstrapping, to explore the finite-sample properties of likelihood-based statistics.    Methodological difficulties, such as the possibile divergence of the likelihood ratio statistic, and prohibitive computational costs limit, however, this possibility to specific applications. 

The review has focused on frequentist hypothesis testing using the likelihood ratio statistic.  Maximum likelihood estimation for a class of nonregular cases, which include the three-parameter Weibull, the gamma, log-gamma and beta distributions, is considered in \cite{Smith85}.  A significant literature has grown since then, parts of which culminated in the book-length account of techniques for parameter estimation in non-standard settings by \cite{Cheng17}.  Most of the difficulties encountered in nonregular settings vanish if the model is analysed using Bayes' rule, though one has always to be cautious.  Bayesian and nonparametric contributions with suitable links to their frequentist counterparts are mentioned in Appendix~B of the Supplementary Material.


\appendix

\section{Prototype demonstrations} 
\label{appendix}

\begin{demo}{Boundary problem \cite[Theorem~3]{SelfLiang87}}
\label{ex:self_liang}
{\rm
Let $y_1, \ldots, y_n$ be $n$ independent observations on the random variable $Y$, and let $l(\theta)$ denote the associated log-likelihood function, where $\theta$ takes values in the parameter space $\Theta$, a subset of ${\mathbb R}^p$.  We want to test whether the true value of $\theta$ lies in the subset of $\Theta$ denoted by $\Theta_0$ versus the alternative that it falls in the complement of $\Theta_0$ in $\Theta$, denoted by $\Theta_1$.  Let $\theta^0$ be the true value of $\theta$, which may fall on the boundary of $\Theta$.  First, expand $2\left\{l(\theta)-l(\theta^0)\right\}$ around $\theta^0$, 
\begin{align*}
2\left\{l(\theta)-l(\theta^0)\right\}&=2(\theta-\theta^0)^\top u(\theta^0)\\&-(\theta-\theta^0)^\top i(\theta^0)(\theta-\theta^0)\\&+O_p(||\theta-\theta^0||^3),
\end{align*}
where $u(\theta)$ is the score function, $i(\theta)$ the Fisher information matrix and $||\cdot||$ represents the Euclidean norm.  Rewrite this expansion as a function of the variable $\tilde Z_n=n^{-1}i_1(\theta^0)^{-1}u(\theta^0)$, where $i(\theta^0)=ni_1(\theta^0)$ and $i_1(\theta^0)$ is the Fisher information matrix associated with a single observation.  This yields
\begin{align*}
2\left\{l(\theta)-l(\theta^0)\right\} = \\
&\hspace{-2cm}- \{\sqrt{n}\tilde Z_n-\sqrt{n}(\theta-\theta^0)\}^\top i_1(\theta^0) \\
& \{\sqrt{n}\tilde Z_n-\sqrt{n}(\theta-\theta^0)\} \\
&\hspace{-2cm}+ u(\theta^0)^\top i(\theta^0)^{-1}u(\theta^0) \\
&\hspace{-2cm}+ O_p(||\theta-\theta^0||^3).
\end{align*}
Consider now the likelihood ratio statistic
\begin{align*}
W & = 2\left\{\sup_{\theta\in\Theta}l(\theta)-\sup_{\theta\in\Theta_0}l(\theta)\right\} \\ 
&\hspace{-1cm} = 
\sup_{\theta\in\Theta}\left[-\{\sqrt{n}\tilde Z_n-\sqrt{n}(\theta-\theta^0)\}^\top i_1(\theta^0)\right. \\
&\left. \{\sqrt{n}\tilde Z_n-\sqrt{n}(\theta-\theta^0)\}\right] \\ 
&\hspace{-1cm} -\sup_{\theta\in\Theta_0}\left[-\{\sqrt{n}\tilde Z_n-\sqrt{n}(\theta-\theta^0)\}^\top i_1(\theta^0)\right. \\
&\left. \{\sqrt{n}\tilde Z_n-\sqrt{n}(\theta-\theta^0)\}\right] \\
&\hspace{-1cm} +O_p(||\theta-\theta^0||^3).
\end{align*}
Approximate the two sets $\Theta$ and $\Theta_0$ by the cones $C_{\Theta-\theta^0}$ and $C_{\Theta_0-\theta^0}$ centered at $\theta^0$,  respectively.  Now, given that $\sqrt{n}\tilde Z_n$ converges in distribution to a multivariate normal distribution with mean zero and covariance matrix $i_1(\theta^0)^{-1}$, for all $\theta$ such that $\theta-\theta^0 = O_p(n^{-1/2})$, the limiting distribution of $W$ becomes
\begin{align*}
\sup_{\theta\in \tilde C}\left\{-(Z-\theta)^\top (Z-\theta)\right\}-\\ \sup_{\theta\in \tilde C_{0}}\left\{-(Z-\theta)^\top (Z-\theta)\right\},
\end{align*}
or equivalently as in Expression~\eqref{selfliangstar}, where $\tilde C$ and $\tilde C_0$ are the corresponding transformations of the cones $C_{\Theta-\theta^0}$ and $C_{\Theta_0-\theta^0}$, respectively, and $Z$ is multivariate standard normal.
}
\end{demo}

\begin{demo}{Non-identifiable parameter \cite[Theorem~2.3]{LiuShao03}}
\label{identifiability_LiuShao}
{\rm
Let $Y_1,\ldots, Y_n$ be $n$ independent and identically distributed random observations from the true distribution function $F^0$. Suppose that we want to test $H_0: \theta \in \Theta_0$ against $H_1:\theta \in\Theta\setminus\Theta_0$, where $\Theta_0=\{\theta\in \Theta : F_{\theta}=F^0\}$ with $F_\theta$ the distribution indexed by $\theta$.  Let
\[
lr(\theta)=\sum_{i=1}^n\log\{\lambda_i(\theta)\}
\] 
be the log-likelihood ratio function, where $\lambda_i(\theta) = \lambda(Y_i;\theta)$ denotes the Radon-Nikodym derivative, $\lambda(\theta)=dF_\theta/dF^0$, evaluated at $Y_i$.  Define the likelihood ratio statistic as
\begin{equation}
\label{W.liushao}
W(H_0) = 2\sup_{\theta\in\Theta\setminus\Theta_0}\{lr(\theta) \vee0 \}, 
\end{equation}
where $\{a\vee b\} = \max(a,b)$.  Assume that there exists a trio $\{S_i(\theta), H(\theta), R_i(\theta)\}$ which satisfies the generalized differentiable in quadratic mean expansion (GDQM)
\[
h_i(\theta)=H(\theta)S_i(\theta)-H^2(\theta)+H^2(\theta)R_i(\theta),
\]
with $h_i(\theta)=\sqrt{\lambda_i(\theta)}-1$.  $H(\theta)$ is the Hellinger distance between $F_\theta$ and $F^0$ defined as
\[ 
H^2(\theta)=E_{F^0}\left[\left\{\sqrt{\lambda_i(\theta)}-1\right\}^2\right]/2
\] 
and $S_i(\theta)$ and $R_i(\theta)$ are such that $E_{F^0}\left[S_i(\theta)\right]=E_{F^0}\left[R_i(\theta)\right]=0$.  Furthermore assume that
\[
\sup_{\theta\in \Theta_{c/\sqrt{n}}}|\nu_n\left(S_i(\theta)\right)|=O_p(1)
\]
and 
\[\sup_{\theta\in \Theta_{c/\sqrt{n}}}|E_{F_n}\left[R_i(\theta)\right]|=o_p(1),\]
for all $c>0$, where $F_n(\cdot)$ indicates the empirical distribution function and $\nu_n(g)=n^{-1/2}(nE_{F_n}-E_{F^0})[g]$ is a random process defined for any integrable function $g$.  Here, $\Theta_{\epsilon}=\{\theta\in\Theta \mid 0<H(\theta)\leq\epsilon\}$ defines the Hellinger neighbourhood of $F^0$.  Now, using the GDQM expansion and a Taylor series expansion of $2\log\{1+h_i(\theta)\}$, the log-likelihood ratio function $lr(\theta)$ can be expressed as
\begin{eqnarray}\label{GDQM_exp}
lr(\theta)&=&2\sum_{i=1}^n\log\{1+h_i(\theta)\}\nonumber\\
&=& 2\sqrt{n}H(\theta)\nu_n(S_i(\theta))\nonumber\\
&-&nH^2(\theta)\left\{2+E_{F_n}\left[S^2_i({\theta})\right]\right\}+o_p(1),
\end{eqnarray}
in $\Theta_{c/\sqrt{n}}$ for all $c > 0$.  Under some general conditions on the trio $\{S_i(\theta),H(\theta),R_i(\theta)\}$ \citep[Theorem~2.2]{LiuShao03}, the quadratic expansion in \eqref{GDQM_exp} holds uniformly in $\theta\in\Theta_{\epsilon}$ for some small enough $\epsilon>0$.  Direct maximization of \eqref{W.liushao} by $\sqrt{n}H(\theta)$ allows us to approximate the likelihood ratio statistic by the quadratic form
\[
\frac{\{\nu_n(S_i(\theta)) \vee 0\}^2}{1+E_{F_n}[S^2_i(\theta)]/2} \approx \{\nu_n(S_i^*(\theta)) \vee 0\}^2.
\]
Let $\mathcal{S}$ be the se of all $\mathcal{L}^2$ limits of the standardized score function 
\[
S_i^*(\theta) = \frac{S_i(\theta)}{\sqrt{1+E_{F^0}[S^2_i(\theta)]/2}}
\]
as $H(\theta)\rightarrow 0$.  To complete the proof we assume there exists a centered Gaussian process $\{G_S:S\in\mathcal{S}\}$ on the same probability space of the empirical process $\nu_n$ with uniformly continuous sample paths and covariance kernel $E_{F^0}\left[G_{S_1}G_{S_2}\right]=E_{F^0}\left[S_1S_2\right]$, for all $S_1,S_2$ belonging to $\mathcal{S}$.  Using results from statistical limit theory, it is possibile to prove the following two  inequalities
\[
W(H_0)\leq \sup_{S\in\mathcal{S}}\{G_S \vee 0\}^2+o_p(1),
\]
\[
W(H_0)\geq \sup_{S\in\mathcal{S}}\{G_S \vee 0\}^2+o_p(1),
\]
which imply that 
\[
\lim_{n\rightarrow \infty}W(H_0)=\sup_{S\in\mathcal{S}}\{G_S \vee 0\}^2.
\]
}
\end{demo}

\begin{demo}{Finite mixture model \cite[Theorem 2.1]{GhoshSen85}}
\label{finite_mixture_GhoshSen85} 
{\rm Let $y_1,\ldots,y_n$ be a sample of $n$ i.i.d.\ observations from the strongly identifiable mixture model~\eqref{2mixture} and
\[
l(\eta)=\sum_{i=1}^n\log\left\{(1-\pi)f_1(y_i;\theta_1)+\pi f_2(y_i;\theta_2)\right\}, 
\]
with $\eta=(\pi,\theta_1,\theta_2)$ be the corresponding log-likelihood function.  Suppose that $H_0:\pi=0$ is true, so the true model density is $f_1(y;\theta_1^0)$, where $\theta_1^0$ is the true value of $\theta_1$.  Unless differently stated, all functions and expectations will be evaluated under this assumption, that is, for $\eta^0=(0,\theta_1^0,\theta_2)$, with arbitrary $\theta_2$.  Let $W(H_0)$ be the likelihood ratio statistic
\begin{eqnarray}
\label{lrmix}
W(H_0) & = & 2\{\sup_{\substack{
	\pi\in [0,1] \\
	\theta_1\in\Theta_1 \\
	\theta_2\in\Theta_2}} l(\eta)-\sup_{\substack{
	\pi=0 \\
	\theta_1\in\Theta_1 \\
	\theta_2\in\Theta_2}}l(\eta)\ \} \nonumber \\
& = &\sup_{\theta_2\in\Theta_2}2\{\sup_{\substack{
	\pi\in [0,1] \\
	\theta_1\in\Theta_1}}l(\eta)-\sup_{\substack{
	\pi=0 \\
	\theta_1\in\Theta_1}}l(\eta) \} .
\end{eqnarray}
Expand $l(\eta)$ with respect to the first two components of $\eta=(\pi, \theta_1, \theta_2)$ around $\pi=0$ and $\theta_1=\theta_1^0$.  This yields
\begin{equation}
\label{GS_expansion}
l(\eta)=l_1(\theta_1^0)+A_n(\eta)+o_p(1),
\end{equation}
where $l_1(\theta_1)=\sum_{i=1}^n\log{f_1(y_i;\theta_1)}$ and 
\begin{eqnarray*}
A_n(\eta) & = &\pi l_{\pi}+(\theta_1-\theta_1^0)^\top l_{\theta_1} + \frac{1}{2}\left\lbrace\pi^2 l_{\pi\pi}\right.\\
			& + &\left. 2\pi (\theta_1-\theta_1^0)^\top l_{\pi\theta_1}\right.\\
			& + &\left.(\theta_1-\theta_1^0)^\top l_{\theta_1\theta_1}(\theta_1-\theta_1^0)\right\rbrace.
\end{eqnarray*}
Here, the two indexes $\pi$ and $\theta_1$ denote differentiation with respect to the corresponding parameter components.  As shown in \cite{GhoshSen85}, in virtue of the Kuhn-Tucker-Lagrange theorem, the unconstrained supremum of $A_n(\eta)$ becomes 
\begin{eqnarray*}
\sup_{\substack{\pi\in [0,1] \\ \theta_1\in\Theta_1}}A_n(\eta) =   
\frac{1}{2} \left\{u_0(\theta_2),u_1^\top\right\} i(\theta_2)^{-1}\left\{u_0(\theta_2),u_1^\top\right\}^\top 
\end{eqnarray*}
if $Z_n(\theta_2) \geq 0$ and 
\[
\sup_{\substack{\pi\in [0,1] \\ \theta_1\in\Theta_1}}A_n(\eta) =\frac{1}{2}u_1^\top i_{11}^{-1}u_1
\]
if $Z_n(\theta_2)<0$, where we define 
$$Z_n(\theta_2)=\frac{\left\{u_{0}(\theta_2)i^{00}(\theta_2)+u_{1}(\theta_2)^\top i^{01}(\theta_2)\right\}}{\{i^{00}(\theta_2)\}^{1/2}}.$$
In the previous three expressions, $u_0(\theta_2)=l_{\pi}(\eta^0)$, $u_1 =l_{\theta_1}(\eta^0)$, $i$ represents the expected information matrix with respect to $\pi$ and $\theta_1$, $i_{jk}(\theta_2)$ denotes the $(jk)$-th component of $i$, for $j=0,1$ and $k=0,1$, while $i^{jk}(\theta_2)$ denotes the $(jk)$-th component of $i^{-1}$.  Similarly, the constrained supremum of $A_n(\eta)$ is 
$$
\sup_{\substack{\pi=0 \\ 
\theta_1\in\Theta_1}}A_n(\eta)=\frac{1}{2}u_1^\top i_{11}^{-1}u_1.
$$
Using known results on the inversion of block matrices, the likelihood ratio statistic \eqref{lrmix} reduces to 
$$
W(H_0) = \sup_{\theta_2\in\Theta_2} Z_n^2(\theta_2)\ I_{\{Z_n\geq 0\}}+ o_{p}(1).
$$
To ensure the convergence of $Z_n(\theta_2)$ to the zero-mean Gaussian processes $Z(\theta_2)$, the set $\Theta_2$ needs be bounded and a Lipschitz condition has to hold for the $u_{0}$ component of the score vector which, in turn, implies tightness of $u_{0}$.  These conditions furthermore guarantee that the remainder term in expansion~(\ref{GS_expansion}) is $o_p(1)$ over the two bounded sets of $\pi$ and $\theta_1$ and uniformly in $\theta_2$.
}
\end{demo}

\begin{demo}{Shift in location for Gaussian model \cite[Theorem 1]{Hawkins77}
\label{gaussian_mean_shift_Hawkins77}}
{\rm 
Given $n$ independent Gaussian observations, we want to test whether
$$Y_i\sim N(\mu,\sigma^2), \quad i=1,\ldots,n,$$
against the alternative that there exists a $0<\tau<n$ at which the unknown mean $\mu$ switches to $\mu^\prime\neq\mu$.  The variance $\sigma^2$ is assumed to be known; we set it to one without loss of generality.  Recall from Section~\ref{shifts-in-location-and-dispersion} that the likelihood ratio statistic can be re-expressed as a function of
\[
U = \max_{1\leq \tau< n}|T_\tau|, 
\]
where
\[
T_\tau=\sqrt{\frac{n}{\tau(n-\tau)}}\sum_{i=1}^\tau(Y_i-\bar{Y}).
\]
The null distribution of $U$ is given at \eqref{u_density}.  The proof considers the following events %
$$
A_\tau=\{|T_\tau|\in(u,u+du)\},
$$
$$
B_\tau=\{|T_i|<|T_\tau|, \forall i\in(1,\ldots,\tau-1)\},
$$
and
$$
C_\tau=\{|T_i|<|T_\tau|, \forall i\in(\tau+1,\ldots,n)\}.
$$
Define
\begin{align*}
& F_U(u+du)-F_U(u)= {\rm Pr}\big\{U\in (u,u+du)\big\}\\
& ={\rm Pr}\left(\bigcup_{\tau=1}^{n-1}\bigg[\{|T_\tau|\in(u,u+du)\}\cap \right. \\ 
& \hspace{20mm} \left. \{|T_\tau|>|T_i|,i\neq \tau\}\bigg]\right)\\
& = \sum_{\tau=1}^{n-1}{\rm Pr}(A_\tau\cap B_\tau\cap C_\tau)\\
& =\sum_{\tau=1}^{n-1}{\rm Pr}(A_\tau){\rm Pr}(B_\tau|A_\tau){\rm Pr}(C_\tau|A_\tau\cap B_\tau).
\end{align*}
Since $T_\tau\sim N(0,1)$, we have that 
$$
{\rm Pr}(A_\tau)=2\phi(u)du+o(du).
$$
Moreover, 
\begin{align}
\nonumber
& {\rm Pr}(B_\tau|A_\tau) \\
\nonumber
& = \ {\rm Pr}(|T_i|<|T_\tau|, \forall i\in(1,\ldots,\tau-1) \mid|T_\tau|=u) \\
\nonumber
& = \ {\rm Pr}(|T_i|<u, \forall i\in(1,\ldots,\tau-1) \mid|T_\tau|=u) +O(du) \\
\nonumber
& = \ g_\tau(u)+O(du), \\
\label{BdA}
\end{align}
\\[-5ex]
where $g_1(u)= 1$ for $u\geq0$ and $g_\tau(u)$ is given in \eqref{eq:g}.
Since the series $\{T_1,T_2,\ldots, T_{n-1}\}$ is Markovian, $\{T_1, T_2, \ldots, T_{\tau-1}\}$ and $\{T_{\tau+1}, T_{\tau+2}, \ldots, T_{n-1}\}$ are independent.  It follows that the events $B_\tau$ and $C_\tau$ are independent given $T_\tau=u$, that is, 
$$
{\rm Pr}(C_\tau|A_\tau \cap B_\tau)={\rm P}(C_\tau|A_\tau).
$$
According to the probability symmetry between $B_\tau$ and $C_\tau$ \citep[\S 2.1.1]{ChenGupta12}, similar to ${\rm Pr}(B_\tau|A_\tau)$, it follows that
\begin{equation}
\label{CdA}
{\rm Pr}(C_\tau|A_\tau)=g_{n-\tau}(u)+O(du).
\end{equation}
Combining \eqref{BdA} and \eqref{CdA}, we obtain
\begin{eqnarray*}
{\rm Pr}\{U\in (u,u+du)\}& = &
2\phi(u)\sum_{\tau=1}^{n-1} g_\tau(u)g_{n-\tau}(u)du \\
    & + & o(du),
\end{eqnarray*}
which corresponds to Expression~\eqref{u_density}.
}
\end{demo}


\section*{Acknowledgements}
Previous and more extended versions of this manuscript can be found on {\tt arXiv:2206.15178}.  We thank the Editors, Associate Editors and all anonymous Referees for their valuable suggestions which greatly helped us improve many aspects of the paper.  It's furthermore a pleasure to acknowledge discussion with Prof. Ruggero Bellio, Prof. Anthony C. Davison and Prof. Nancy Reid.  This research was supported by University of Padova grant no. CPDA101912 \textit{Large- and small-sample inference under non-standard conditions} (Progetto di Ricerca di Ateneo 2010).


\renewcommand\refname{References}


\clearpage

\renewcommand\refname{{\small SUPPLEMENTARY MATERIAL} \\[3ex] B~~~Annotated bibliography}

\renewcommand{\bibpreamble}{{
\footnotesize The subsequent section-wise list of references supplements the work cited in the main text in an attempt to provide a comprehensive overview of the asymptotic properties of the likelihood ratio statistic in nonregular problems. \\ ~}}

\end{document}